\documentclass{gtart_h}  

\def\ifplaintex{\expandafter\ifx\csname documentclass\endcsname\relax}

\ifplaintex 
\else
\fi

\expandafter\ifx\csname beginpicture\endcsname\relax
\expandafter\ifx\csname documentclass\endcsname\relax
\input pictex \else
\input prepictex \input pictex \input postpictex \fi\fi

\def\gt{{\mathsurround=0pt\it $\cal G\mskip-2mu$eometry \&\ 
$\cal T\!\!$opology}}        

\def\gtp{{\mathsurround=0pt\it $\cal G\mskip-2mu$eometry \&\ 
$\cal T\!\!$opology $\cal P\!$ublications}}  

\def\lognumber#1{\def\thelognumber{#1}}
\def\volumenumber#1{\def\thevolumenumber{#1}}
\def\papernumber#1{\def\thepapernumber{#1}}
\def\volumeyear#1{\def\thevolumeyear{#1}}

\def\pagenumbers#1#2{\def\startpage{#1}\def\finishpage{#2}}
\def\published#1{\def\publishdate{#1}}
\def\proposed#1{\def\theproposer{#1}}
\def\seconded#1{\def\theseconders{#1}}
\def\received#1{\def\receiveddate{#1}}

\def\accepted#1{\def\accepteddate{#1}}
\def\asciititle#1{\def\theasciititle{#1}}
\def\covertitle#1{\def\thecovertitle{#1}}
\def\coverauthors#1{\def\thecoverauthors{#1}}
\def\asciiauthors#1{\def\theasciiauthors{#1}}
\def\asciiaddress#1{\def\theasciiaddress{#1}}

\long\def\asciiabstract#1{\long\def\theasciiabstract{#1}}
\def\asciikeywords#1{\def\theasciikeywords{#1}}

\def\shorttitle#1{\def\theshorttitle{#1}}

\let\\\par\let\thelognumber\relax
\let\thevolumenumber\relax\let\thepapernumber\relax
\let\thevolumeyear\relax\let\thesamplenumber\relax\let\startpage\relax
\let\finishpage\relax\let\publishdate\relax\let\receiveddate\relax
\let\reviseddate\relax\let\accepteddate\relax\let\theasciititle\relax
\let\thecovertitle\relax\let\theasciiauthors\relax\let\theasciiaddress\relax
\let\theasciiabstract\relax\let\theasciikeywords\relax
\let\theasciiemail\relax\let\theshortauthors\relax\let\theshorttitle\relax
\let\thecoverauthors\relax

\long\def\maketitlep{   

\count0=\startpage

\gt\hfill      
\beginpicture
\setcoordinatesystem units <0.33truein, 0.33truein> point at 2.2 0.9
\setplotsymbol ({$\cal G$})
\plotsymbolspacing=9truept
\circulararc 315 degrees from 0 1 center at 0 0
\setplotsymbol ({$\cal T$})
\circulararc 315 degrees from 1 -1 center at 1 0
\endpicture
\break
{\small\ifx\thesamplenumber\relax 
Volume \else Sample
\fi\thevolumenumber\ (\thevolumeyear)
\startpage--\finishpage\nl
Published: \publishdate}
\vglue 0.5truein plus 0.4fil minus 0.1truein

{\parskip=0pt\leftskip 0pt plus 1fil\def\\{\par\smallskip}{\ifplaintex\large
\else\Large\fi\bf\thetitle}\par\medskip}   

\vglue 0pt plus 0.1fil 

{\parskip=0pt\leftskip 0pt plus 1fil\def\\{\par}{\sc\theauthors}
\par\medskip}

\vglue 0pt plus 0.1fil 

{\small\parskip=0pt\let\newline\\
{\leftskip 0pt plus 1fil\def\\{\par}{\sl\theaddress}\par}
\expandafter\ifx\theemail\relax    
\relax\else\vglue 5pt plus 0.02fil minus 2pt\def\\{\stdspace{\rm 
and}\stdspace} 
\cl{Email:\stdspace\tt\theemail}\fi
\ifx\theurl\relax                  
\relax\else\vglue 5pt plus 0.02fil minus 2pt\def\\{\stdspace{\rm 
and}\stdspace}
\cl{URL:\stdspace\tt\theurl}\fi\par}

\vglue 7pt plus 0.3fil minus 3pt

{\bf Abstract}
\vglue 5pt plus 0.1fil minus 2pt

\theabstract

\vglue 7pt plus 0.3fil minus 3pt

{\bf AMS Classification numbers}\quad Primary:\quad \theprimaryclass

Secondary:\quad \thesecondaryclass

\vglue 5pt plus 0.3fil minus 2pt

{\bf Keywords:}\quad \thekeywords

\vglue 10pt plus 0.5fil minus 5pt

{\small  Proposed: \theproposer\hfill Received: \receiveddate\nl
Seconded: \theseconders\hfill 
\ifx\reviseddate\relax                         
Accepted: \accepteddate                        
\else
Revised: \reviseddate                          
\fi}
\eject
}       

\let\maketitlepage\maketitlep
\let\maketitle\maketitlepage

\font\phead=cmsl9 scaled 950
\font=cmsl9 scaled 1050
\font\pnum=cmbx10 scaled 913
\font\lnum=cmbx10 
\font\pfoot=cmsl9 scaled 950
\font=cmsl9 scaled 1050
\ifplaintex
\headline{\vbox to 0pt{\vskip -4.5mm\line{\small\phead\ifnum
\count0=\startpage ISSN 1364-0380 (on line)
1465-3060 (printed) \hfill {\pnum\folio}\else\ifodd\count0\def\\{ }%
\ifx\theshorttitle\relax\thetitle\else\theshorttitle\fi\hfill{\pnum\folio}
\else\def\\{ and }{\pnum\folio}\hfill\ifx\theshortauthors\relax\theauthors
\else\theshortauthors\fi\fi\fi}\vss}}
\footline{\vbox to 0pt{\vglue 0mm\line{\small\pfoot\ifnum\count0=\startpage
\copyright\ \gtp\hfill\else
\gt, Volume \thevolumenumber\ (\thevolumeyear)\hfill\fi}\vss
}}
\else
\makeatletter
\def\@oddhead{{\small\lhead\ifnum\count0=\startpage ISSN 1364-0380 (on line)
1465-3060 (printed) \hfill {\lnum\number\count0}\else\ifodd\count0
\def\\{ }\ifx\theshorttitle\relax \thetitle \else\theshorttitle\fi\hfill
{\lnum\number\count0}\else\def\\{ and }{\lnum\number\count0}
\hfill\ifx\theshortauthors\relax 
\theauthors\else\theshortauthors\fi\fi\fi}}\def\@evenhead{\@oddhead}
\def\@oddfoot{\small\lfoot\ifnum\count0=\startpage\copyright\ \gtp\hfill\else
\gt, Volume \thevolumenumber\ (\thevolumeyear)\hfill\fi}
\def\@evenfoot{\@oddfoot}
\makeatother
\fi

\newwrite\gtoutfile
\long\gdef\makeheadfile{  
{\def\\{, }\def\s{ }
\immediate\openout\gtoutfile head.xxx
\immediate\write\gtoutfile{Proxy-for: \ifx\theasciiauthors\relax
\theauthors\else\theasciiauthors\fi\s<\ifx\theasciiemail\relax\theemail\else\theasciiemail\fi>}
\immediate\write\gtoutfile{\noexpand\\}
\immediate\write\gtoutfile{Authors: \ifx\theasciiauthors\relax
\theauthors\else\theasciiauthors\fi}
{\def\\{ }\immediate\write\gtoutfile{Title: \ifx\theasciititle\relax
\thetitle\else\theasciititle\fi}}
\immediate\write\gtoutfile{Subj-class: GT or SG or MG etc}
\immediate\write\gtoutfile{MSC-class: \theprimaryclass\ifx\thesecondaryclass\relax\else, \thesecondaryclass\fi}
\immediate\write\gtoutfile{Journal-ref: Geom. Topol. \thevolumenumber
(\thevolumeyear) \startpage-\finishpage}
\immediate\write\gtoutfile{Comments: Published by Geometry and Topology at}
\immediate\write\gtoutfile{\s\s http://www.maths.warwick.ac.uk/gt/GTVol\thevolumenumber/paper\thepapernumber.abs.html}
\immediate\write\gtoutfile{\noexpand\\}
\immediate\write\gtoutfile{}
\ifx\theasciiabstract\relax
\immediate\write\gtoutfile{\theabstract}\else
\immediate\write\gtoutfile{\theasciiabstract}\fi
\immediate\write\gtoutfile{}
\immediate\write\gtoutfile{\noexpand\\}
\immediate\write\gtoutfile{}
\immediate\closeout\gtoutfile}}  

\def\maketitlepage{\maketitlep\makeheadfile}
\let\maketitle\maketitlepage

\lognumber{539}
\received{8 December 2004}
\volumenumber{9}\papernumber{37}\volumeyear{2005}
\pagenumbers{1639}{1676}   
\published{28 August 2005}
\accepted{19 August 2005}
\proposed{Gunnar Carlsson}
\seconded{Ralph Cohen, Bill Dwyer}

\usepackage{amsmath,amsfonts,amscd,amssymb}

\allowdisplaybreaks

\DeclareMathAlphabet\EuR{U}{eur}{m}{n}
\SetMathAlphabet\EuR{bold}{U}{eur}{b}{n}

\newcommand{\calfin}{{\mathcal F}\!{\mathcal I}\!{\mathcal N}}

\newcommand{\calicof}{{\mathcal I}\!{\mathcal C}\!{\mathcal O}\!{\mathcal F}}

\newcommand{\calvcyc}{{\mathcal V}\!{\mathcal C}\!{\mathcal Y}\!{\mathcal C}}
\newcommand{\caltr}{{\mathcal T}\!{\mathcal R}}

\newcommand{\calf}{{\cal F}}
\newcommand{\calg}{{\cal G}}
\newcommand{\calh}{{\cal H}}

\newcommand{\bbC}{{\mathbb C}}

\newcommand{\bbP}{{\mathbb P}}
\newcommand{\bbQ}{{\mathbb Q}}
\newcommand{\bbR}{{\mathbb R}}

\newcommand{\bbZ}{{\mathbb Z}}

\newcommand{\bfK}{{\mathbf K}}
\newcommand{\bfL}{{\mathbf L}}

\newcommand{\asmb}{\operatorname{asmb}}
\newcommand{\aut}{\operatorname{aut}}

\newcommand{\Hei}{\operatorname{Hei}}
\newcommand{\id}{\operatorname{id}}
\newcommand{\ind}{\operatorname{ind}}

\newcommand{\pr}{\operatorname{pr}}

\newcommand{\res}{\operatorname{res}}

\newcommand{\tors}{\operatorname{tors}}
\newcommand{\UNil}{\operatorname{UNil}}
\newcommand{\Wh}{\operatorname{Wh}}

\newcommand{\incl}{\operatorname{incl}}

\newcommand{\pt}{\{\ast\}}

\newtheorem{theorem}{Theorem}[section]
\newtheorem{lemma}[theorem]{Lemma}

\newtheorem{remark}[theorem]{Remark}
\newtheorem{corollary}[theorem]{Corollary}

{\catcode`@=11\global\let\c@equation=\c@theorem}

\renewcommand{\theenumi}{\rm(\roman{enumi})}

\newcommand{\EGF}[2]{E_{#2}(#1)}               

\newcommand{\comsquare}[8]                   
{\begin{CD}
#1 @>#2>> #3\\
@V{#4}VV @VV{#5}V\\
#6 @>>#7> #8
\end{CD}
}

\begin{document}

\title{$K$-- and $L$--theory of the semi-direct product of 
the\\discrete 3--dimensional Heisenberg group by $\bbZ/4$}
\shorttitle{$K$-- and $L$--theory of $\Hei \rtimes \bbZ/4$}
\asciititle{K- and L-theory of the semi-direct product of the
discrete 3-dimensional Heisenberg group by Z/4}
\covertitle{$K$-- and $L$--theory of the semi-direct product\\of the discrete 
3--dimensional Heisenberg group by ${\noexpand \bf Z}/4$}

\author{Wolfgang L\"uck}
\asciiauthors{Wolfgang Lueck}
\coverauthors{Wolfgang L\noexpand\"uck}
\address{Fachbereich Mathematik, Universit\"at 
M\"unster\\Einsteinstr.~62, 48149 M\"unster, Germany}
\asciiaddress{Fachbereich Mathematik, Universitaet 
Muenster\\Einsteinstr. 62, 48149 Muenster, Germany}
\email{lueck@math.uni-muenster.de}
\urladdr{www.math.uni-muenster.de/u/lueck/}

\begin{abstract}
We compute the group homology, the topological $K$--theory of the
reduced $C^*$--algebra, the algebraic $K$--theory and the algebraic
$L$--theory of the group ring of the semi-direct product of the
three-dimensional discrete Heisenberg group by $\bbZ/4$.  These
computations will follow from the more general treatment of a certain
class of groups $G$ which occur as extensions $1 \to K \to G \to Q \to
1$ of a torsionfree group $K$ by a group $Q$ which satisfies certain
assumptions. The key ingredients are the Baum--Connes and Farrell--Jones
Conjectures and methods from equivariant algebraic topology.
\end{abstract}
\asciiabstract{%
We compute the group homology, the topological K-theory of the reduced
C^*-algebra, the algebraic K-theory and the algebraic L-theory of the
group ring of the semi-direct product of the three-dimensional
discrete Heisenberg group by Z/4.  These computations will follow from
the more general treatment of a certain class of groups G which occur
as extensions 1-->K-->G-->Q-->1 of a torsionfree group K by a group
Q which satisfies certain assumptions. The key ingredients are the
Baum-Connes and Farrell-Jones Conjectures and methods from
equivariant algebraic topology.}

\asciikeywords{K- and L-groups of group rings and group C^*-algebras,
three-dimensional Heisenberg group.}

\maketitle

\setcounter{section}{-1}
\section{Introduction}

The original motivation for this paper was the question of Chris Phillips
how the topological $K$--theory of the reduced (complex) $C^*$--algebra of the semi-direct product
$\Hei \rtimes \bbZ/4$ looks like.
Here $\Hei$ is the \emph{three-dimensional discrete
Heisenberg group} which is the subgroup of $GL_3(\bbZ)$ consisting of
upper triangular matrices with $1$ on the diagonals. The $\bbZ/4$--action is given by:
$$\left(\begin{array}{ccc} 1& x & y \\ 0 & 1 & z \\ 0 & 0 & 1\end{array}\right)
\mapsto
\left(\begin{array}{ccc} 1& -z & y-xz \\ 0 & 1 & x \\ 0 & 0 & 1\end{array}\right)$$
The answer, which is proved in Theorem~\ref{topological $K$--theory of Hei rtimes Z/4},
consists of an explicit isomorphism
$$
j_0 \bigoplus c[0]_0' \bigoplus c[2]_0' \colon
K_0(\pt) \bigoplus \widetilde{R}_{\bbC}(\bbZ/4) \bigoplus \widetilde{R}_{\bbC}(\bbZ/2)
\xrightarrow{\cong}  K_0(C^*_r(\Hei \rtimes \bbZ/4))
$$
and a short exact sequence
$$
0 \to \widetilde{R}_{\bbC}(\bbZ/4) \bigoplus \widetilde{R}_{\bbC}(\bbZ/2)
\xrightarrow{c[0]_1' \bigoplus c[2]_1'} K_1(C^*_r(\Hei \rtimes \bbZ/4)) \xrightarrow{c_1}
\widetilde{K}_1(S^3) \to 0,
$$
which splits since $\widetilde{K}_1(S^3) \cong \bbZ$.
Here $\widetilde{R}_{\bbC}(\bbZ/m)$ is the kernel of the split surjective map
$R_{\bbC}(\bbZ/m) \to R_{\bbC}(\{1\}) \cong \bbZ$ which sends the class of a complex $\bbZ/m$--representation
to the class of $\bbC \otimes_{\bbC[\bbZ/m]} V$.
As abelian group we get for $n \in \bbZ$
$$K_n(C^*_r(\Hei \rtimes \bbZ/4))  \cong \bbZ^5.$$
This computation will play a role in the paper
by Echterhoff, L\"uck and Phillips~\cite{Echterhoff-Lueck-Phillips(2005)},
where certain $C^*$--algebras given by semi-direct products of rotation algebras
with finite cyclic groups are classified.

Although the group $\Hei \rtimes \bbZ/4$ is
very explicit, this computation is highly non-trivial
and requires besides the Baum--Connes Conjecture a lot of machinery from
equivariant algebraic topology. Even harder is the computation of the
middle and lower $K$--theory. The result is
(see Corollary~\ref{cor: algebraic K-theory of G = Hei rtimes Z/4})
\begin{eqnarray*}
\Wh_n(\Hei \rtimes \bbZ/4) & \cong & \left\{\begin{array}{lll}
N\!K_n(\bbZ [\bbZ/4]) \bigoplus N\!K_n(\bbZ [\bbZ/4])  & & \text { for  } n = 0,1;
\\
0 & & \text{ for } n \le -1,
\end{array}\right.
\end{eqnarray*}
where $N\!K_n(\bbZ [\bbZ/4])$ denotes the $n$-th Nil-group of
$\bbZ[\bbZ/4]$ which appears in the Bass--Heller--Swan decomposition of
$\bbZ[\bbZ/4 \times \bbZ]$.
So the lower $K$--theory is trivial and the middle $K$--theory is completely
made up of Nil-groups.

We also treat the $L$--groups. The answer and calculation
is rather messy due to the appearance  of
$\UNil$--terms and the structure of the family of infinite virtually cyclic subgroups
(see Theorem~\ref{the: L-theory of hei rtimes Z/4}). If one is willing to invert $2$,
these $\UNil$--terms and questions about decorations disappear and the answer is
given by the short split exact sequence:
\begin{multline*}
0 \to L_n(\bbZ)\left[\frac{1}{2}\right] \bigoplus
\widetilde{L}_n(\bbZ[\bbZ/2])\left[\frac{1}{2}\right] \bigoplus
\widetilde{L}_{n-1}(\bbZ[\bbZ/2])\left[\frac{1}{2}\right]
\\ \bigoplus \widetilde{L}_n(\bbZ[\bbZ/4])\left[\frac{1}{2}\right] \bigoplus
\widetilde{L}_{n-1}(\bbZ[\bbZ/4])\left[\frac{1}{2}\right]
\\
\xrightarrow{j} L_n(\bbZ[\Hei \rtimes \bbZ/4])\left[\frac{1}{2}\right]
\to L_{n-3}(\bbZ)\left[\frac{1}{2}\right] \to 0
\end{multline*}
Finally we will also compute the group homology
(see Theorem~\ref{the: group homology})
\begin{eqnarray*}
H_n(G) & = & \bbZ/2 \times \bbZ/4 \text{ for } n \ge 1, n \not= 2,3;
\\
H_2(G) & = & \bbZ/2;
\\
H_3(G) & = & \bbZ \times \bbZ/2 \times \bbZ/4.
\end{eqnarray*}
In turns out that we can handle a much more general setting provided
that the Baum--Connes Conjecture or the Farell--Jones Conjecture is true for $G$. Namely, we
will consider an extension of (discrete) groups
\begin{eqnarray}
& 1 \to K \xrightarrow{i}  G \xrightarrow{p} Q \to 1 &
\label{extension of groups}
\end{eqnarray}
which satisfies the following conditions:
\begin{description}

\item[\rm(M)] Each non-trivial finite subgroup of $Q$ is contained
in a unique maximal finite subgroup;

\item[\rm(NM)] Let $M$ be a maximal finite subgroup of $Q$. Then $N_QM = M$ unless
$G$ is torsionfree;

\item[\rm(T)] $K$ is torsionfree.

\end{description}

The special case, where $K$ is trivial, is treated in \cite[Theorem 5.1]{Davis-Lueck(2003)}.
In \cite[page 101]{Davis-Lueck(2003)} it is explained using
\cite[Lemma 4.5]{Lueck-Stamm(2000)}),
\cite[Lemma 6.3]{Lueck-Stamm(2000)} and
\cite[Propositions 5.17, 5.18 and 5.19 in II.5 on pages 107 and 108]{Lyndon-Schupp(1977)}
why the following groups satisfy conditions (M) and (NM):
\begin{itemize}

\item Extensions $1 \to \bbZ^n \to Q \to F \to 1$ for finite $F$ such that the conjugation
  action of $F$ on $\bbZ^n$ is free outside $0 \in \bbZ^n$;

\item Fuchsian groups;

\item One-relator groups.

\end{itemize}
Of course $\Hei \rtimes \bbZ/4$ is an example for $G$. For such groups $G$ we will
establish certain exact Mayer--Vietoris sequences relating the $K$-- or $L$--theory
of $G$ to the $K$-- and $L$--theory of $p^{-1}(M)$ for maximal finite subgroups $M \subseteq Q$
and terms involving the quotients $G\backslash \underline{E}G$ and
$p^{-1}(M)\backslash \underline{E}p^{-1}(M)$. The classifying space $\underline{E}G$
for proper $G$--actions plays an important role
and often there are nice small geometric models for them.
One key ingredient in the computations for
$\Hei \rtimes \bbZ/4$ will be to show that
$G\backslash \underline{E}G$ in this case is $S^3$. For instance the computation
of the group homology illustrates that it is often very convenient
to work with the spaces $G\backslash \underline{E}G$ although one wants information about
$BG$.

\section{Topological $K$--theory}
\label{sec: Topological K-theory of the reduced C^*-algebras of certain classes of groups}

For a $G$--$CW$--complex $X$ let $K^G_*(X)$ be its \emph {equivariant $K$--homology theory}.
If $G$ is trivial, we abbreviate $K_*(X)$. For a $C^*$--algebra $A$ let $K_*(A)$ be its
topological $K$--theory. Recall that a model $\underline{E}G$ for the
\emph{classifying space for proper $G$--actions} is  a $G$--$CW$--complex with finite
isotropy groups such that $(\underline{E}G)^H$ is contractible for each finite subgroup $H
\subseteq G$. It has the property that for any $G$--$CW$--complex $X$ with finite isotropy
groups there is precisely one $G$--map from $X$ to $\underline{E}G$ up to $G$--homotopy.
In particular two models for $\underline{E}G$ are $G$--homotopy equivalent.
For more information about the spaces $\underline{E}G$ we refer for instance
to \cite{Baum-Connes-Higson(1994)}, \cite{Lueck(2004h)}, \cite{Mislin(2003)}.
\cite{Dieck(1987)}.
Recall that the \emph{Baum--Connes Conjecture}
(see \cite[Conjecture 3.15 on page 254]{Baum-Connes-Higson(1994)}) says that
the assembly map
$$\asmb\colon K_n^G(\underline{E}G) \xrightarrow{\cong} K_n(C^*_r(G))$$
is an isomorphism for each $n \in \bbZ$, where $C^*_r(G)$ is the
\emph{reduced group $C^*$--algebra} associated to $G$. (For an identification of the assembly map used
in this paper with the original one we refer to
Hambleton--Pedersen~\cite{Hambleton-Pedersen(2004)}). Let $EG$ be a model for the
\emph{classifying space for free $G$--actions}, ie,
a free $G$--$CW$--complex which is contractible (after forgetting the group action).
Up to $G$--homotopy there is precisely one $G$--map
$s \colon EG \to \underline{E}G$. The \emph{classical assembly map}
$a$ is defined as the composition
$$a \colon K_n(BG) = K_n^G(EG) \xrightarrow{K_p^G(s)}
K_n^G(\underline{E}G) \xrightarrow{\asmb} K_n(C_r^*(G)).$$
For more information about the Baum--Connes Conjecture we refer for instance
to \cite{Baum-Connes-Higson(1994)}, \cite{Lueck-Reich(2003b)}, \cite{Mislin(2003)},
\cite{Valette(2002)}.

From now on consider a group $G$ as described in \eqref{extension of groups}
We want to compute $K_n^G(\underline{E}G)$. If $G$ satisfies the Baum--Connes Conjecture
this is the same as $K_n(C^*_r(G))$.

First we construct a nice model for $\underline{E}Q$. Let
$\{(M_i) \mid i \in I\}$ be the set of conjugacy classes of maximal
finite subgroups of $M_i \subseteq Q$. By attaching free $Q$--cells we get
an inclusion of $Q$--$CW$--complexes
$j_1 \colon \coprod_{i \in I} Q \times_{M_i} EM_i  \to EQ$.
Define $\underline{E}Q$ as the $Q$--pushout
\begin{eqnarray}
& \comsquare{\coprod_{i \in I} Q \times_{M_i} EM_i}{j_1}{EQ}
{u_1}{f_1}{\coprod_{i \in I} Q/M_i}{k_1}{\underline{E}Q} &
\label{pushout for underline E Q}
\end{eqnarray}
where $u_1$ is the obvious $Q$--map obtained by collapsing each $EM_i$ to a point.

We have to explain why $\underline{E}Q$ is a model for the classifying space for
proper actions of $Q$. Obviously it is a $Q$--$CW$--complex. Its  isotropy groups
are all finite. We have to show for $H \subseteq Q$ finite
that $(\underline{E}Q)^H$ contractible. We begin with the case $H \not= \{1\}$.
Because of conditions (M) and (NM)
there is precisely one index $i_0 \in I$ such that $H$ is subconjugated to
$M_{i_0}$ and is not subconjugated to $M_i$ for $i \not= i_0$ and we get
$$\left(\coprod_{i \in I} Q/M_i\right)^H~=~\left(Q/M_{i_0}\right)^H~=~\pt.$$
It remains to treat $H = \{1\}$. Since $u_1$ is a non-equivariant homotopy equivalence
and $j_1$ is a cofibration, $f_1$ is a non-equivariant homotopy equivalence and hence
$\underline{E}Q$ is contractible (after forgetting the group action).

Let $X$ be a $Q$--$CW$--complex and $Y$ be a $G$--$CW$--complex. Then
$X \times Y$ with the $G$--action given by $g \cdot (x,y) = (p(g)x,gy)$ is a
$G$--$CW$--complex and the $G$--isotropy group $G_{(x,y)}$ of $(x,y)$ is
$p^{-1}(H_x) \cap G_y$. Hence $\underline{E}Q \times \underline{E}G$ is a
$G$--$CW$--model for $\underline{E}G$ and $EQ \times \underline{E}G$ is a $G$--$CW$--model
for $EG$, since $\ker(p \colon G \to Q)$ is torsionfree by assumption.
Let $Z$ be a $M_i$--$CW$--complex. Then there is a $G$--homeomorphism
$$G \times_{p^{-1}(M_i)} \left(Z \times  \res_G^{p^{-1}(M_i)} Y\right) \xrightarrow{\cong}
(Q \times_{M_i} Z) \times Y
\hspace{5mm} (g,(z,y)) \mapsto ((p(g),z),gy).$$
The inverse sends $((q,z),y)$ to $(g,(z,g^{-1}y)$ for any choice of $g \in G$ with $p(g) = q$.
If we cross the $Q$--pushout \eqref{pushout for underline E Q} with $\underline{E}G$,
then we obtain the following $G$--pushout:
\begin{eqnarray}
& \comsquare{\coprod_{i \in I} G \times_{p^{-1}(M_i)} Ep^{-1}(M_i)}{j_2}{EG}
{u_2}{f_2}{\coprod_{i \in I} G \times_{p^{-1}(M_i)} \underline{E}p^{-1}(M_i)}{k_2}{\underline{E}G}
\label{pushout for underline E G}
\end{eqnarray}
If we divide out the $G$--action in the pushout
\eqref{pushout for underline E G} above we obtain the pushout:
\begin{eqnarray}
& \comsquare{\coprod_{i \in I} Bp^{-1}(M_i)}{j_3}{BG}
{u_3}{f_3}{\coprod_{i \in I} p^{-1}(M_i)
\backslash\underline{E}p^{-1}(M_i)}{k_3}{G\backslash\underline{E}G}
\label{pushout for G backslash underline E G}
\end{eqnarray}
If we divide out the $Q$--action in the pushout
\eqref{pushout for underline E Q} we obtain the pushout:
\begin{eqnarray}
& \comsquare{\coprod_{i \in I} BM_i}{j_4}{BQ}
{u_4}{f_4}{\coprod_{i \in I} \pt}{k_4}{Q\backslash\underline{E}Q} &
\label{pushout for Q backslash underline E Q}
\end{eqnarray}

\begin{theorem} \label{the: Mayer-Vietoris for topological K-theory}
Let $G$ be the group appearing in \eqref{extension of groups} and assume
that conditions (M), (NM) and  (T) hold. Assume that $G$ and all
groups $p^{-1}(M_i)$ satisfy the Baum--Connes Conjecture. Then the
Mayer--Vietoris sequence associated to
\eqref{pushout for underline E G}  yields the long exact sequence of abelian groups:
\begin{multline*}
\ldots \xrightarrow{\partial_{n+1}}
\bigoplus_{i \in I} K_n(Bp^{-1}(M_i))
\\
\xrightarrow{\left(\bigoplus_{i \in I} K_n(Bl_i)\right)
\bigoplus \left(\bigoplus_{i \in I}  a[i]_n \right) }
 K_n(BG)  \bigoplus \left(\bigoplus_{i \in I} K_n(C_r^*(p^{-1}(M_i))\right)
\\
\xrightarrow{a_n  \bigoplus \left(\bigoplus_{i \in I}  K_n(C_r^*(l_i))\right)}
K_n(C^*_r(G)) \xrightarrow{\partial_n}
\bigoplus_{i \in I} K_{n-1}(Bp^{-1}(M_i))
\\
\xrightarrow{\left(\bigoplus_{i \in I} K_{n-1}(Bl_i)\right)\bigoplus
\left(\bigoplus_{i \in I}  a[i]_{n-1} \right) }
 K_{n-1}(BG) \bigoplus  \left(\bigoplus_{i \in I} K_{n-1}(C_r^*(p^{-1}(M_i)) \right)
\\
\xrightarrow{a_{n-1}  \bigoplus \left(\bigoplus_{i \in I}  K_{n-1}(C_r^*(l_i))\right)} \ldots
\end{multline*}
Here the maps $a[i]_n$ and $a$ are classical assembly
maps and  $l_i \colon p^{-1}(M_i) \to G$ is the
inclusion.

Let $\Lambda$ be a ring with $\bbZ \subseteq \Lambda \subseteq\bbQ$ such that the order
of each finite subgroup of $G$ is invertible in $\Lambda$. Then the
composition
\begin{multline*}
\Lambda \otimes_{\bbZ} K_n(Bp^{-1}(M_i))
\xrightarrow{\id_{\Lambda} \otimes_{\bbZ}  a[i]_n}
K_n(C_r^*(p^{-1}(M_i)) = K_n^{p^{-1}(M_i)}(\underline{E}p^{-1}(M_i))
\\
\xrightarrow{\id_{\Lambda} \otimes_{\bbZ}  \ind_{p^{-1}(M_i) \to \{1\}}}
\Lambda \otimes_{\bbZ} K_n(p^{-1}(M_i)\backslash \underline{E}p^{-1}(M_i))
\end{multline*}
is an isomorphism, where $\ind$ denotes the induction map.
In particular the long exact sequence above reduces
after applying $\Lambda \otimes_{\bbZ} $ to split exact
short exact sequences of  $\Lambda$--modules:
\begin{multline*}
0 \to
\bigoplus_{i \in I} \Lambda \otimes_{\bbZ}  K_n(Bp^{-1}(M_i))
\xrightarrow{\left(\bigoplus_{i \in I} \id_{\Lambda} \otimes_{\bbZ} K_n(Bl_i)\right)
\bigoplus \left(\bigoplus_{i \in I}  \id_{\Lambda} \otimes_{\bbZ}  a[i]_n \right)}
\\
\Lambda \otimes_{\bbZ}  K_n(BG) \bigoplus
 \left(\bigoplus_{i \in I} \Lambda \otimes_{\bbZ}  K_n(C_r^*(p^{-1}(M_i))) \right)
\\
\xrightarrow{\bigoplus_{i \in I}  \id_{\Lambda} \otimes_{\bbZ}
K_n(C_r^*(l_i))  \bigoplus \id_{\Lambda} \otimes_{\bbZ}  a_n}
\Lambda \otimes_{\bbZ}  K_n(C^*_r(G)) \to 0
\end{multline*}
\end{theorem}
\begin{proof} The Mayer Vietoris sequence is obvious using the fact that for
a free $G$--$CW$--complex $X$ there is a canonical isomorphism
$K_n^G(X) \xrightarrow{\cong} K_n(G\backslash X)$.
The composition
\begin{multline*}
\Lambda \otimes_{\bbZ} K_n(Bp^{-1}(M_i))
\xrightarrow{\id_{\Lambda} \otimes_{\bbZ}  a[i]_n}
K_n(C_r^*(p^{-1}(M_i)) = K_n^{p^{-1}(M_i)}(\underline{E}p^{-1}(M_i))
\\
\xrightarrow{\id_{\Lambda} \otimes_{\bbZ}  \ind_{p^{-1}(M_i) \to \{1\}}}
\Lambda \otimes_{\bbZ} K_n(p^{-1}(M_i)\backslash \underline{E}p^{-1}(M_i))
\end{multline*}
is bijective by \cite[Lemma 2.8 (a)]{Lueck-Stamm(2000)}.
\end{proof}

The advantage of the following version is that it involves the spaces
$G\backslash \underline{E}G$ instead of the spaces $BG$, and these often have rather small geometric models.
In the case $G = \Hei\rtimes \bbZ/4$ we will see that $G\backslash \underline{E}G$
is the three-dimensional sphere $S^3$ (see Lemma~\ref{G backslash underline E G is S^3}).

\begin{theorem}
\label{the: modified Mayer-Vietoris sequence}
Let $G$ be the group appearing in \eqref{extension of groups} and assume
conditions (M), (NM) and (T) hold. Assume that $G$ and all
groups $p^{-1}(M_i)$ satisfy the Baum--Connes Conjecture. Then
there is a long exact sequence of abelian groups:
\begin{multline*}
\ldots \xrightarrow{c_{n+1} \bigoplus_{i \in I} d[i]_{n+1}}
K_{n+1}(G\backslash\underline{E}G) \xrightarrow{\partial_{n+1}}
\bigoplus_{i \in I} K_n(C_r^*(p^{-1}(M_i)))
\\
\xrightarrow{\left(\bigoplus_{i \in I} K_n(C_r^*(l_i))\right) \bigoplus
\left(\bigoplus_{i \in I} c[i]_n\right)}
K_n(C_r^*(G)) \bigoplus
\left(\bigoplus_{i \in I} K_{n}(p^{-1}(M_i)\backslash\underline{E}p^{-1}(M_i))\right)
\\
\xrightarrow{c_{n} \bigoplus_{i \in I} d[i]_{n}}
K_{n}(G\backslash\underline{E}G) \xrightarrow{\partial_{n}}
\bigoplus_{i \in I} K_{n-1}(C_r^*(p^{-1}(M_i)))
\\
\xrightarrow{\left(\bigoplus_{i \in I} K_{n-1}(C_r^*(l_i)) \right)\bigoplus
\left(\bigoplus_{i \in I} c[i]_{n-1}\right)} \ldots
\end{multline*}
Here the homomorphisms $d[i]_n$ come from the 
$\left(p^{-1}(M_i) \to G\right)$--equivariant maps $\underline{E}p^{-1}(M_i) \to \underline{E}G$
which are unique up to equivariant homotopy.
The maps $c_n$ and (analogously for $c[i]_n$) are the compositions
$$K_n(C^*_r(G)) \xrightarrow{\asmb^{-1}} K_n^G(\underline{E}G) \xrightarrow{\ind_{G \to \{1\}}}
K_n(G\backslash \underline{E}G).$$
Let $\Lambda$ be a ring with $\bbZ \subseteq \Lambda \subseteq\bbQ$ such that the order
of each finite subgroup of $G$ is invertible in $\Lambda$. Then the
composition
\begin{multline*}
\Lambda \otimes_{\bbZ} K_n(BG)
\xrightarrow{\id_{\Lambda} \otimes_{\bbZ} a_n} \Lambda \otimes_{\bbZ}  K_n(C^*_r(G))
\xrightarrow{\id_{\Lambda} \otimes_{\bbZ}  c_{n}}
\Lambda \otimes_{\bbZ} K_n(G\backslash \underline{E}G)
\end{multline*}
is an isomorphism of $\Lambda$--modules. In particular the long exact sequence above
reduces after applying $\Lambda \otimes_{\bbZ} -$ to split exact short sequences of
$\Lambda$--modules:
\begin{multline*}
0 \to \bigoplus_{i \in I} \Lambda \otimes_{\bbZ} K_n(C_r^*(p^{-1}(M_i)))
\xrightarrow{\left(\bigoplus_{i \in I}
\id_{\Lambda} \otimes_{\bbZ} K_n(C_r^*(l_i))\right) \bigoplus
\left(\bigoplus_{i \in I} \id_{\Lambda} \otimes_{\bbZ} c[i]_n\right)}
\\
\Lambda \otimes_{\bbZ} K_n(C_r^*(G)) \bigoplus \left(\bigoplus_{i \in I}
\Lambda \otimes_{\bbZ} K_{n}(p^{-1}(M_i)\backslash\underline{E}p^{-1}(M_i))\right)
\\
\xrightarrow{\id_{\Lambda} \otimes_{\bbZ} c_{n}
\bigoplus_{i \in I} \id_{\Lambda} \otimes_{\bbZ} d[i]_{n}}
\Lambda \otimes_{\bbZ} K_{n}(G\backslash\underline{E}G)  \to 0
\end{multline*}
\end{theorem}
\begin{proof}
From the pushout \eqref{pushout for G backslash underline E G} we get the long exact
Mayer Vietoris sequence for (non-equivariant) topological $K$--theory
\begin{multline*}
\ldots \xrightarrow{\partial_{n+1}}
\bigoplus_{i \in I} K_n(Bp^{-1}(M_i))
\xrightarrow{\left(\bigoplus_{i \in I} K_n(Bl_i)\right) \bigoplus
\left(\bigoplus_{i \in I} K_n((p^{-1}(M_i)\backslash s_i)\right)}
\\
K_n(BG) \bigoplus \left(\bigoplus_{i \in I}
K_n(p^{-1}(M_i)\backslash \underline{E}p^{-1}(M_i))\right)
\xrightarrow{H_n(G\backslash s) \bigoplus \left(\bigoplus_{i \in I} d[i]_n\right)}
K_n(G\backslash\underline{E}G)
\\
\xrightarrow{\partial_n}
\bigoplus_{i \in I} K_{n-1}(Bp^{-1}(M_i))
\xrightarrow{\left(\bigoplus_{i \in I} K_{n-1}(Bl_i)\right) \bigoplus
\left(\bigoplus_{i \in I} K_{n-1}(p^{-1}(M_i)\backslash s_i)\right)}
\\
K_{n-1}(BG) \bigoplus \left(\bigoplus_{i \in I}
K_{n-1}(p^{-1}(M_i)\backslash \underline{E}p^{-1}(M_i))\right)
\xrightarrow{K_{n-1}(G\backslash s) \bigoplus \left(\bigoplus_{i \in I} d[i]_{n-1}\right)} \ldots
\end{multline*}
where $s_i \colon Ep^{-1}(M_i) \to \underline{E}p^{-1}(M_i)$ and
$s \colon EG \to \underline{E}G$ are  (up to equivariant homotopy unique) equivariant maps.
Now one splices the long exact Mayer--Vietoris sequences  from above and from
Theorem~\ref{the: Mayer-Vietoris for topological K-theory} together.
\end{proof}

\section{The semi-direct product of the Heisenberg group and a cyclic group of order four}
\label{sec: The example Hei rtimes bbZ/4}

We want to study the following example.
Let $\Hei$ be the \emph{discrete Heisenberg group}. We will use the presentation
\begin{eqnarray}
\Hei &  = &\langle u,v,z \mid [u,z] = 1, [v,z] = 1, [u,v] = z\rangle.
\label{presentation of Hei}
\end{eqnarray}
Throughout this section let $G$ be the semi-direct product
$$G = \Hei \rtimes \bbZ/4$$
with respect to the homomorphism $\bbZ/4 \to \aut(\Hei)$
which sends the generator $t$ of $\bbZ/4$
to the automorphism of $\Hei$ given on generators by
$z \mapsto z$, $u\mapsto v$ and $v \mapsto u^{-1}$.
Let $Q$ be the semi-direct product
$\bbZ^2 \rtimes \bbZ/4$ with respect to the automorphism
$\bbZ^2 \to \bbZ^2$ which comes from multiplication with the complex number
$i$ and the inclusion $\bbZ^2 \subseteq \bbC$. Since the action of
$\bbZ/4$ on $\bbZ^2$ is free outside $0$, the group $Q$ satisfies (M) and (NM)
(see \cite[Lemma 6.3]{Lueck-Stamm(2000)}).
The group  $G$ has the presentation
$$G  = \langle u,v,z,t \mid [u,z] = [v,z] = [t,z] = t^4 = 1,
[u,v] = z, tut^{-1} = v, tvt^{-1} = u^{-1}\rangle.$$
Let $i \colon \bbZ \to G$ be the inclusion sending the generator of
$\bbZ$ to $z$. Let $p \colon G \to Q$ be the group homomorphism,
which sends $z$ to the unit element, $u$ to $(1,0)$ in $\bbZ^2 \subseteq Q$,
$v$ to $(0,1)$ in $\bbZ^2 \subseteq Q$ and $t$ to the
generator of $\bbZ/4 \subseteq Q$. Then $1 \to \bbZ \to G \to Q \to 1$ is a
central extension which satisfies the conditions (M), (NM) and (T)
appearing in \eqref{extension of groups}.
Moreover, $G$ is amenable and hence $G$ and all its subgroups satisfy the Baum--Connes
Conjecture \cite{Higson-Kasparov(2001)}.

In order to apply the general results above we have to figure out
the conjugacy classes of finite subgroups of
$Q = \bbZ^2 \rtimes \bbZ/4$ and among them the maximal ones.
An element of order $2$ in $Q$ must  have the form $xt^2$ for $x \in \bbZ^2$.
In the sequel we write the group multiplication in $Q$ and $G$
multiplicatively and in $\bbZ^2$ additively.
We compute $(xt^2)^2 = xt^2xt^2 = (x-x) = 0$. Hence the set of elements of order two in
$Q$ is $\{xt^2\mid x \in \bbZ^2\}$. Consider $e_1 = (1,0)$ and $e_2 = (0,1)$ in $\bbZ^2$.
We claim that up to conjugacy
there are the following subgroups of order two:
$\langle e_1t^2 \rangle, \langle e_1e_2t^2\rangle, \langle t^2\rangle$.
This follows from the computations for $x,y \in \bbZ^2$
\begin{eqnarray*}
y(xt^2)y^{-1} & = & yxyt^2 = (x+2y)t^2;
\\
t(xt^2)t^{-1} & = & txt^{-1}t^2 = (ix)t^2.
\end{eqnarray*}
An element of order $4$ must have the form $xt$ for $x \in \bbZ^2$. We compute
$$(xt)^4 = xtxt^{-1}t^2xt^{-2}t^3xt^{-3} = (x +ix +i^2x + i^3x) = (1 +i +i^2+i^3)x = 0x = 0.$$
Hence the set of elements of order four in
$Q$ is $\{xt\mid x \in \bbZ^2\}$.
We claim that up to conjugacy
there are the following subgroups of order four:
$\langle e_1t \rangle, \langle t\rangle$.
This follows from the computations for $x,y \in \bbZ^2$
\begin{eqnarray*}
y(xt)y^{-1} & = & (x + y -iy)t;
\\
t(xt)t^{-1} & = & (ix)t.
\end{eqnarray*}
We have $(e_1t)^2 = e_1te_1t = e_1ie_1t^2 = e_1 e_2t^2$.
The considerations above imply:

\begin{lemma} \label{lem: maximal finite subgroups of Q}
Up to conjugacy $Q$ has the following non-trivial finite subgroups
$$\langle e_1t^2 \rangle, \langle e_1e_2t^2\rangle, \langle t^2\rangle,
\langle e_1t \rangle, \langle t\rangle.$$
The maximal finite subgroups are up to conjugacy
$$M_0 = \langle t\rangle , M_1 = \langle e_1t \rangle, M_2 = \langle e_1t^2 \rangle.$$
\end{lemma}

Since $t^4 = 1$, $(ut^2)^2 = ut^2ut^{-2} = uu^{-1} = 1$ and
$(ut)^4 = utut^{-1}t^2ut^{-2}t^3ut^{-3} = uvu^{-1}v^{-1} = z$ hold in $G$,
the preimages of these groups under $p \colon G \to Q$ are given by
\begin{eqnarray*}
p^{-1}(M_0) & = & \langle t,z \rangle  \cong  \bbZ/4 \times \bbZ;
\\
p^{-1}(M_1) & = & \langle ut, z \rangle  = \langle ut \rangle  \cong  \bbZ;
\\
p^{-1}(M_2) & = & \langle ut^2,z \rangle  \cong  \bbZ/2 \times \bbZ.
\end{eqnarray*}
One easily checks

\begin{lemma} \label{lem: finite subgroups of G}
Up to conjugacy the finite subgroups of $G$ are
$\langle t \rangle$, $\langle t^2 \rangle$ and
$\langle ut^2 \rangle$.
\end{lemma}

Next we construct nice geometric models for $\underline{E}G$ and its
orbit space $G\backslash \underline{E}G$. Let $\Hei(\bbR)$ be the \emph{real Heisenberg group},
ie, the Lie group of real $(3,3)$--matrices of the special form:
$$\left(\begin{array}{ccc} 1& x & y \\ 0 & 1 & z \\ 0 & 0 & 1\end{array}\right)$$
In the sequel we identify such a matrix with the
element $(x,y,z) \in \bbR^3$. Thus $\Hei(\bbR)$ can be identified with the Lie group
whose underlying manifold is $\bbR^3$ and whose group multiplication is given by
$$(a,b,c) \bullet (x,y,z) = (a+x,b + y +az, c+z).$$
The discrete Heisenberg group is given by the subgroup where all the
entries $x,y,z$ are integers. In the presentation of the discrete Heisenberg group
\eqref{presentation of Hei} the elements $u$, $v$ and $z$  correspond to $(1,0,0)$,
$(0,0,1)$ and $(0,1,0)$. Obviously $\Hei$ is
a torsionfree discrete subgroup of the contractible Lie
group $\Hei(\bbR)$. Hence $\Hei(\bbR)$ is a model for $E\!\Hei$ and $\Hei\backslash \Hei(\bbR)$
for $B\!\Hei$. We have the following $\bbZ/4$--action on $\Hei(\bbR)$,
with the generator $t$ acting
by $(x,y,z) \mapsto (-z,y-xz,x)$. This is an action by automorphisms of Lie groups and
induces the homomorphism $\bbZ/4 \to \aut(\Hei)$ on $\Hei$ which we have used above to define
$G = \Hei \rtimes \bbZ/4$. The $\Hei$--action and $\bbZ/4$--action on $\Hei(\bbR)$ above fit
together to a $G = \Hei \rtimes \bbZ/4$--action. The next result is the main geometric
input for the desired computations.

\begin{lemma}
\label{G backslash underline E G is S^3}
The manifold $\Hei(\bbR)$ with the $G$--action above is a model for $\underline{E}G$.
The quotient space $G\backslash \underline{E}G$ is homeomorphic to $S^3$.
\end{lemma}
\begin{proof}
Let $\bbR \subseteq \Hei(\bbR)$ be
the subgroup of elements $\{(0,y,0) \mid y \in \bbR\}$. This is the center of $\Hei(\bbR)$.
The intersection $\bbR \cap \Hei$ is $\bbZ \subseteq \bbR$. Thus we get a $\bbR/\bbZ =
S^1$--action on $\Hei\backslash \Hei(\bbR)$. One easily checks that this
$S^1$--action and the $\bbZ/4$--action above commute so that we see
a $S^1 \times \bbZ/4$--action on $\Hei\backslash \Hei(\bbR)$. The $S^1$--action is free,
but the $S^1 \times \bbZ/4$--action is not. Next we figure out its fixed points.

Obviously $t^2$ sends $(x,y,z)$ to
$(-x,y,-z)$. We compute for $(a,b,c) \in \Hei$, $u \in \bbR$ and
$(x,y,z) \in \Hei(\bbR)$
\begin{eqnarray*}
(a,b,c) \cdot (0,u,0) \cdot t \cdot (x,y,z) & = & (a-z,u +b+y-xz-ax,c+x);\\
(a,b,c) \cdot (0,u,0) \cdot t^2 \cdot (x,y,z) & = &  (a-x,u + b+y-az,c-z); \\
(a,b,c) \cdot (0,u,0) \cdot (x,y,z) & = & (a+x,u + b+y+az,c+z).
\end{eqnarray*}
Hence the isotropy group of $\Hei \cdot (x,y,z) \in
\Hei\backslash\Hei(\bbR)$ under the $S^1 \times \bbZ/4$--action
contains $(\exp(2\pi i u),t)$ in its isotropy
group under the $S^1 \times \bbZ/4$--action if and only if
$(a-z,u + b+y-xz-ax,c+x) = (x,y,z)$ holds for some integers $a,b,c$.
The last statement is equivalent to the condition
that $2x$ and $x + z$  are integers, $y$ is an arbitrary real number and
$u - 3 x^2 \in \bbZ$.

The isotropy group of $\Hei \cdot (x,y,z) \in
\Hei\backslash\Hei(\bbR)$ contains $(\exp(2\pi i u),t^2)$ in its isotropy
group under the $S^1 \times \bbZ/4$--action if and only if
$(a-x,u + b+y-az,c-z) = (x,y,z)$ holds for some integers $a,b,c$.
Obviously the last statement is equivalent to the condition
that $2x$, $2z$ and $u- 2xz$ are integers and $y$ is an arbitrary real number.

The isotropy group of $\Hei \cdot (x,y,z) \in
\Hei\backslash\Hei(\bbR)$ contains $(\exp(2\pi i u),1)$ in its isotropy
group under the $S^1 \times \bbZ/4$--action if and only if
$(a+x,u + b+y+az,c+z) = (x,y,z)$ holds for some integers $a,b,c$.
The last statement is equivalent to the condition
that $x = 0$, $z = 0$, $u$ is an integer and $y$ is an arbitrary real number.

This implies that the orbits under the $S^1 \times \bbZ/4$--action on
$\Hei\backslash\Hei(\bbR^3)$ are free except the orbits through
$\Hei\cdot (1/2,0,1/2)$, whose isotropy group is
the cyclic subgroup of order four generated by $(\exp(3\pi i/4),t)$,
and the orbits though $\Hei \cdot (0,0,0)$, whose isotropy group is
the cyclic subgroup of order four generated by $(\exp(0),t)$,
and the orbits though $\Hei \cdot (1/2,0,0)$ and $\Hei \cdot (0,0,1/2)$, whose isotropy
groups are the cyclic subgroup of order two generated by $(\exp(0),t^2)$. By the
slice theorem any point $p \in \Hei\backslash\Hei(\bbR)$ has a
neighborhood of the form $S^1 \times \bbZ/4 \times_{H_p} U_p$, where
$H_p$ is its isotropy group and $U_p$ a $2$--dimensional real
$H_p$--representation, namely the tangent space of $\Hei\backslash\Hei(\bbR)$ at $p$.
Since there are only finitely $S^1 \times \bbZ/4$--orbits
which are non-free, the $H_p$--action on $U_p$ is free outside the origin for each
$p \in \Hei\backslash \Hei(\bbR)$. In particular $H_p \backslash U_p$ is a manifold
without boundary.
If the isotropy group $H_p$ is mapped under the
projection $\pr \colon S^1 \times \bbZ/4\to S^1$ to the trivial group, then
$\bbZ/4\backslash\left(S^1 \times \bbZ/4 \times_{H_p} U_p\right)$ is $S^1$--homeomorphic
to $S^1 \times H_p\backslash U_p$ and hence a free $S^1$--manifold without boundary.
If the projection $\pr \colon S^1 \times \bbZ/4 \to S^1$  is injective on $H_p$,
then $\bbZ/4\backslash\left(S^1 \times \bbZ/4 \times_{H_p} U_p\right)$ is the
$S^1$--manifold $S^1 \times_{H_p} U_p$ with respect to the free $H$--action on $S^1$ induced by $p$
which has no boundary and precisely one non-free $S^1$--orbit.
This shows that the quotient of $\Hei\backslash\Hei(\bbR^3)$ under the $\bbZ/4$--action
is a closed $S^1$--manifold with precisely one non-free orbit.

The fixed point set of any finite subgroup of $G$ of the $G$--space
$\Hei(\bbR) = \bbR^3$ is a non-empty affine real subspace of
$\Hei(\bbR) = \bbR^3$ and hence contractible. This shows that
$\Hei(\bbR)$ with its  $G$--action is a model for $\underline{E}G$.
Hence $G\backslash \underline{E}G$ is a  closed $S^1$--manifold
with precisely one non-free orbit, whose quotient space under
the $S^1$--action is the orbit space of $T^2$ under the
$\bbZ/4$--action. One easily checks for the rational homology
$$H_n\left((\bbZ/4)\backslash T^2;\bbQ\right) \cong_{\bbQ}
H_n(T^2)\otimes_{\bbZ[\bbZ/4]} \bbQ \cong H_n(S^2;\bbQ).$$
This implies that the $S^1$--space $G\backslash\underline{E}G$ is a Seifert bundle over
$(\bbZ/4)\backslash T^2 \cong S^2$ with precisely one singular fiber.
Since the orbifold fundamental group of this orbifold $S^2$ with precisely one
cone point vanishes, the map $e \colon \pi_1(S^1) \to \pi_1(G\backslash\underline{E}G)$ given
by evaluating the $S^1$--action at some base point is surjective by
\cite[Lemma 3.2] {Scott(1983)}.
The Hurewicz map $h \colon \pi_1(G\backslash \underline{E}G) \to
H_1(G\backslash \underline{E}G) $ is
bijective since $\pi_1(G\backslash\underline{E}G)$ is a quotient
of $\pi_1(S^1)$ and hence is abelian. The composition
$$\pi_1(S^1) \xrightarrow{e} \pi_1(G\backslash\underline{E}G)
\xrightarrow{h} H_1(G\backslash \underline{E}G)$$
agrees with the composition
$$\pi_1(S^1) \xrightarrow{h'} H_1(S^1) =  H_1(\bbZ\backslash \bbR) \xrightarrow{e'}
H_1(\Hei\backslash \Hei(\bbR)) \xrightarrow{H_1(\pr)} H_1(G\backslash
\underline{E}G),$$
where $h'$ is the Hurewicz map, $e'$ given by evaluating the
$S^1$--operation and $\pr$ is the obvious projection.
The map $H_1(\bbZ\backslash \bbR) \to H_1(\Hei\backslash \Hei(\bbR))$ is trivial
since the element $z \in \Hei$ is a commutator, namely $[u,v]$.
Hence $G\backslash \underline{E}G$ is a simply connected closed Seifert fibered
$3$--manifold. We conclude from \cite[Lemma 3.1]{Scott(1983)}
that $G\backslash\underline{E}G$ is homeomorphic to $S^3$.
\end{proof}

Next we investigate what information
Theorem~\ref{the: modified Mayer-Vietoris sequence} gives in combination with
Lemma~\ref{G backslash underline E G is S^3}.

We have to analyze the maps
$$c[i]_n \colon K_n(C_r^*(p^{-1}(M_i)))  \to
K_{n}(p^{-1}(M_i)\backslash\underline{E}p^{-1}(M_i)),$$
which are defined as the compositions
\begin{multline*}
K_n(C^*_r(p^{-1}(M_i))) \xrightarrow{\asmb^{-1}}
K_n^{p^{-1}(M_i)}(\underline{E}p^{-1}(M_i))
\\
\xrightarrow{\ind_{p^{-1}(M_i) \to \{1\}}}
K_n(p^{-1}(M_i)\backslash \underline{E}p^{-1}(M_i)).
\end{multline*}
For $i = 1$ the group $p^{-1}(M_i)$ is isomorphic to $\bbZ$ and hence the maps
$c[1]_n$ are all isomorphisms. In the case $i = 0,2$ the group
$p^{-1}(M_i)$ looks like $H_i \times \bbZ$ for $H_0 = \langle t \rangle \cong \bbZ/4$ and
$H_2 = \langle ut^2 \rangle \cong \bbZ/2$. The following diagram commutes:
$$\begin{CD}
K_n(C_r^*(H_i \times \bbZ)) @< \asmb < \cong < K_n^{H_i \times \bbZ}(\underline{E}H_i \times \bbZ)
@< \cong << K_n^{H_i}(\pt) \bigoplus K_{n-1}^{H_i}(\pt)
\\
@VK_n(C_r^*(\pr_i) )VV @V\ind_{H_i \times \bbZ \to \bbZ}VV
@V \ind_{H_i \to \{1\}}  \bigoplus \ind_{H_i \to \{1\}}VV
\\
K_n(C_r^*(\bbZ)) @< \asmb < \cong < K_n^{\bbZ}(E\bbZ)
@< \cong << K_n(\pt) \bigoplus K_{n-1}(\pt)
\end{CD}$$
The map $\ind_{H_i \to \{1\}} \colon K_n^{H_i}(\pt) \to K_n(\pt)$ is the map
$\id \colon 0 \to 0$ for $n$ odd. For $n$ even it can be identified with the  homomorphism
$\epsilon\colon R_{\bbC}(H_i) \to \bbZ$ which sends the class of a
complex $H_i$--representation $V$ to the complex dimension of $\bbC \otimes_{\bbC H_i} V$.
This map is split surjective.
The kernel of $\epsilon$ is denoted by $\widetilde{R}_{\bbC}(H_i)$.
Define for $i = 0,2$ maps
\begin{eqnarray}
c[i]_n' \colon \widetilde{R}_{\bbC}(H_i) & \to & K_n(C^*_r(G))
\label{the maps c[i]_n'}
\end{eqnarray}
as follows. For $n$ even it is the composition
$$\widetilde{R}_{\bbC}(H_i) \subseteq R_{\bbC}(H_i)  =
K_n(C_r^*(H_i)) \xrightarrow{K_n(C_r^*(l_i'))}
K_n(C_r^*(G)),$$
where $l_i' \colon H_i \to G$ is the inclusion.
For $n$ odd it is the composition
$$\widetilde{R}_{\bbC}(H_i) \subseteq R_{\bbC}(H_i)  = K_{n-1}(C_r^*(H_i))
\xrightarrow{x_i} K_n(C^*_r(H_i\times \bbZ))  \xrightarrow{K_n(C_r^*(l_i))}
K_n(C_r^*(G)),$$
where $l_i \colon H_i \times \bbZ = p^{-1}(M_i) \to G$ is the inclusion and
$$x_i \bigoplus  K_n(y_i) \colon K_{n-1}(C_r^*(H_i))
\bigoplus K_n(C_r^*(H_i)) \xrightarrow{\cong}
K_n(C^*_r(H_i \times \bbZ))$$
is the canonical isomorphism for $y_i \colon H_i \to H_i \times \bbZ$ the inclusion.
The map
$$\partial_n \colon K_{n}(G\backslash\underline{E}G) \to
\bigoplus_{i \in I} K_{n-1}(C_r^*(p^{-1}(M_i))) $$
appearing in Theorem~\ref{the: modified Mayer-Vietoris sequence}
vanishes after applying $\bbQ \otimes_{\bbZ} -$. Since the  target
is a finitely generated torsionfree abelian group, the map itself is trivial.
Hence we obtain from Theorem~\ref{the: modified Mayer-Vietoris sequence}
short exact sequences for $n \in \bbZ$
\begin{multline*}
0 \to \widetilde{R}_{\bbC}(\bbZ/4) \bigoplus \widetilde{R}_{\bbC}(\bbZ/2)
\xrightarrow{c[0]_n' \bigoplus c[2]_n'} K_n(C^*_r(G)) \xrightarrow{c_n} K_n(S^3) \to 0,
\end{multline*}
where we identify $H_0 = \langle t \rangle = \bbZ/4$ and
$H_2 = \langle ut^2 \rangle = \bbZ/2$ and $G\backslash\underline{E}G = S^3$  using
Lemma~\ref{G backslash underline E G is S^3}. If
$j_n \colon K_n(\pt) = K_n(C^*_r(\{1\})) \to K_n(C^*_r(G))$
is induced by the inclusion of the trivial subgroup,
we can rewrite the sequence above as the short
exact sequence
\begin{multline*}
0 \to K_n(\pt) \bigoplus \widetilde{R}_{\bbC}(\bbZ/4)
\bigoplus \widetilde{R}_{\bbC}(\bbZ/2)
\xrightarrow{j_n \bigoplus c[0]_n' \bigoplus c[2]_n'} K_n(C^*_r(G)) \xrightarrow{c_n}
\\
\widetilde{K}_n(S^3) \to 0,
\end{multline*}
where $\widetilde{K}_n(Y)$ is for a path connected space $Y$
the cokernel of the obvious map $K_n(\pt) \to K_n(Y)$.
We have $\widetilde{K}_0(S^3) = 0$ and $\widetilde{K}_1(S^3) \cong \bbZ$.
Thus we get
\begin{theorem}
\label{topological $K$--theory of Hei rtimes Z/4}
We have the isomorphism
$$
j_0 \bigoplus c[0]_0' \bigoplus c[2]_0' \colon
K_0(\pt) \bigoplus \widetilde{R}_{\bbC}(\bbZ/4) \bigoplus \widetilde{R}_{\bbC}(\bbZ/2)
 \xrightarrow{\cong}  K_0(C^*_r(\Hei \rtimes \bbZ/4))
$$
and the short exact sequence
$$
0 \to \widetilde{R}_{\bbC}(\bbZ/4) \bigoplus \widetilde{R}_{\bbC}(\bbZ/2)
\xrightarrow{c[0]_1' \bigoplus c[2]_1'} K_1(C^*_r(G)) \xrightarrow{c_1}
\widetilde{K}_1(S^3) \to 0,
$$
where the maps $c[i]_n'$ have been defined in \eqref{the maps c[i]_n'}.
In particular $K_n(C_r^*(G))$ is a free abelian group of rank five for all $n$.
\end{theorem}

\begin{remark} \label{rem: consistent with Chern characters} \em
These computations are consistent with
the computation of $K_n(C^*_r(G))\left[\frac{1}{2}\right]$
coming from the Chern character constructed in \cite{Lueck(2002d)}.
\em\end{remark}

\begin{remark} \label{real $K$--theory} \em
One can also use these methods to compute the topological $K$--theory
of the real reduced group $C^*$--algebra $C_r^*(\Hei \rtimes \bbZ/4;\bbR)$.
One obtains the short exact sequence
\begin{multline*}
0 \to KO_n(\pt) \bigoplus
\widetilde{K}_n(C_r^*(\bbZ/2;\bbR)) \bigoplus \widetilde{K}_{n-1}(C_r^*(\bbZ/2;\bbR))
\\
\bigoplus
\widetilde{K}_n(C_r^*(\bbZ/4;\bbR)) \bigoplus \widetilde{K}_{n-1}(C_r^*(\bbZ/4;\bbR))
\\
\to K_n(C_r^*(\Hei \rtimes \bbZ/4)) \to \widetilde{KO}_n(S^3) \to 0,
\end{multline*}
which splits after inverting $2$.
\em
\end{remark}

\section{Algebraic $K$--theory}
\label{sec: Algebraic $K$-theory}

In this section we want to describe
what the methods above yield for the algebraic $K$--theory
provided that instead of the Baum--Connes
Conjecture the relevant version of the Farrell--Jones Conjecture for
algebraic $K$--theory (see \cite{Farrell-Jones(1993a)}) is true.
The $L$--theory will be treated in the next section.
We want to prove the following:

\begin{theorem} \label{lem: result for algebraic $K$-theory}
Let $R$ be a regular ring, for instance $R = \bbZ$.
Let $G$ be the group appearing in \eqref{extension of groups} and assume
that conditions (M), (NM), and (T)  are satisfied.
Suppose that $G$ and all subgroups $p^{-1}(M_i)$ satisfy the Farrell--Jones Conjecture
for algebraic $K$--theory with coefficients in $R$. Then we get for $n \in \bbZ$  the isomorphism
$$\bigoplus_{i \in I} \Wh_n(Rl_i) \colon \Wh_n(R[p^{-1}(M_i)])
\xrightarrow{\cong} \Wh_n(RG),$$
where $l_i \colon p^{-1}(M_i) \to G$ is the inclusion.
\end{theorem}

Notice that in the context of the Farrell--Jones Conjecture
one has to consider the \emph{family of virtually cyclic subgroups} $\calvcyc$ and only
under special assumptions it suffices to consider the \emph{family
$\calfin$ of finite subgroups}. Recall that a family $\calf$ of subgroups
is a set of subgroups closed under conjugation and taking subgroups
and that a model for \emph{the classifying space $\EGF{G}{\calf}$ for the family
$\calf$} is a $G$--$CW$--complex whose isotropy groups
belong to $\calf$ and whose $H$--fixed point set
is contractible for each $H \in \calf$. It is characterized up to
$G$--homotopy by the property that any $G$--$CW$--complex, whose isotropy
groups belong to $\calf$, possesses up to $G$--homotopy precisely one
$G$--map to $\EGF{G}{\calf}$. In particular two models for
$\EGF{G}{\calf}$ are $G$--homotopy equivalent and for an inclusion of
families $\calf \subset\calg$ there is up to $G$--homotopy precisely
one $G$--map $\EGF{G}{\calf} \to \EGF{G}{\calg}$. The space
$\underline{E}G$ is the same as $\EGF{G}{\calfin}$.

Let $\calh^G_*(X;\bfK(R?))$ and $\calh^G_*(X;\bfL^{\langle -\infty \rangle}(R?))$ be the
$G$--homology theories associated to the algebraic $K$ and $L$--theory
spectra over the orbit category
$\bfK(R?)$ and $\bfL^{\langle -\infty \rangle}(R?)$ (see \cite{Davis-Lueck(1998)}).
They satisfy for each subgroup $H \subseteq G$
\begin{eqnarray*}
\calh^G_n(G/H;\bfK(R?)) & \cong & K_n(RH);
\\
\calh^G_n(G/H;\bfL^{\langle -\infty \rangle}(R?)) & \cong &
L_n^{\langle -\infty \rangle}(RH).
\end{eqnarray*}
The \emph{Farrell--Jones Conjecture} (see \cite[1.6 on page 257]{Farrell-Jones(1993a)})
says that the projection $\EGF{G}{\calvcyc} \to G/G$ induces
isomorphisms
\begin{eqnarray*}
\calh^G_n(\EGF{G}{\calvcyc};\bfK(R?)) & \xrightarrow{\cong} &
\calh^G_n(G/G;\bfK(R?)) = K_n(RG);
\\
\calh^G_n(\EGF{G}{\calvcyc};\bfL^{\langle - \infty \rangle}(R?)) & \xrightarrow{\cong} &
\calh^G_n(G/G;\bfL^{\langle - \infty \rangle}(R?)) = L^{\langle - \infty \rangle}_n(RG).
\end{eqnarray*}
In the $L$--theory case one must use $L^{\langle -\infty\rangle}$.
There are counterexamples to the Farrell--Jones Conjecture for the other decorations
$p$, $h$ and $s$ (see \cite{Farrell-Jones-Lueck(2002)}).

In the sequel we denote for a $G$--map $f \colon X \to Y$ by
$\calh^G_n(f \colon X \to Y;\bfK(R?))$ the value of $\calh^G_n$ on the
pair given by the mapping cylinder of $f$ and $Y$ viewed as a
$G$--subspace. We will often use the long exact sequence associated to
this pair:
\begin{multline*}
\ldots \to \calh^G_n(X;\bfK(R?)) \to \calh^G_n(Y;\bfK(R?)) \to
\calh^G_n(f \colon X \to Y;\bfK(R?))
\\
\to \calh^G_{n-1}(X;\bfK(R?)) \to \calh^G_{n-1}(Y;\bfK(R?)) \to \ldots
\end{multline*}
The following result is taken from \cite{Bartels(2003h)}.
\begin{theorem} \label{the: Bartels results}
There are isomorphisms
\begin{multline*}
\calh^G_n(\underline{E}G;\bfK(\bbR?)) \bigoplus
\calh^G_n\left(\underline{E}G \to \EGF{G}{\calvcyc};\bfK(R?)\right)
\\
\xrightarrow{\cong}  \calh^G_n(\EGF{G}{\calvcyc};\bfK(R?));
\end{multline*}
\begin{multline*}
\calh^G_n(\underline{E}G;\bfL^{\langle -\infty \rangle}(\bbR?)) \bigoplus
\calh^G_n\left(\underline{E}G \to \EGF{G}{\calvcyc};\bfL^{\langle -\infty \rangle}(R?)\right)
\\
\xrightarrow{\cong}  \calh^G_n(\EGF{G}{\calvcyc};\bfL^{\langle - \infty  \rangle}(R?)),
\end{multline*}
where in the $K$--theory context $G$ and $R$ are arbitrary and in the $L$--theory context
$G$ is arbitrary and we assume for any virtually cyclic subgroup $V \subseteq G$ that
$K_{-i}(RV) = 0$ for sufficiently large $i$.
\end{theorem}

For a virtually cyclic group $V$ we have
$K_{-i}(\bbZ V) = 0$ for $n \ge 2$ (see \cite{Farrell-Jones(1995)}).

The terms $\calh^G_n\left(\underline{E}G \to \EGF{G}{\calvcyc};\bfK(R?)\right)$
vanish for instance if $R$ is a regular ring containing $\bbQ$.
The terms $\calh^G_n\left(\underline{E}G\to
\EGF{G}{\calvcyc};\bfL^{\langle - \infty \rangle}(R?)\right)$
vanish after inverting $2$ (see Lemma~\ref{lem: general simplification of calvcyc}).
 Recall that the Whitehead group $\Wh_n(RG)$ by definition is
 $\calh^G_n(EG \to G/G;\bfK(R?))$. This implies that
$\Wh_n(RG) = \calh^G_n(EG \to \EGF{G}{\calvcyc};\bfK(R?))$ if
the Farrell--Jones Isomorphism Conjecture for algebraic $K$--theory holds
for $RG$. The group $\Wh_1(\bbZ G)$ is the classical Whitehead group $\Wh(G)$.
If $R$ is a principal ideal domain, then
$\Wh_0(RG)$ is $\widetilde{K}_0(RG)$ and
$\Wh_n(RG) = K_n(RG)$ for $n \le -1$.

If we cross the $Q$--pushout \eqref{pushout for underline E Q}
with $\EGF{G}{\calvcyc}$ we obtain the $G$--pushout:
\begin{eqnarray}
&
\comsquare{\coprod_{i \in I} G \times_{p^{-1}(M_i)}
\EGF{p^{-1}(M_i)}{\calvcyc(K \cap p^{-1}(M_i))}}
{}{\EGF{G}{\calvcyc(K)}}
{}{}
{\coprod_{i \in I} G \times_{p^{-1}(M_i)}
\EGF{p^{-1}(M_i)}{\calvcyc}}
{}{\EGF{G}{\calvcyc_f}}
&
\label{G-pushout for EGcalvcyc_f}
\end{eqnarray}
where $\calvcyc(K \cap p^{-1}(M_i))$ is the family of virtually cyclic subgroups of
$p^{-1}(M_i)$, which are contained in $K \cap p^{-1}(M_i)$,
and $\calvcyc(K)$ is the family of virtually cyclic subgroups of
$G$, which are contained in $K$, and $\calvcyc_f$ is the family of virtually cyclic
subgroups of $G$, whose image under $p\colon G \to Q$ is finite.
Since $K$ is torsionfree,
elements in $\calvcyc(K \cap p^{-1}(M_i))$ and
$\calvcyc(K)$ are trivial or infinite cyclic groups.
The following result is taken from \cite[Theorem 2.3]{Lueck-Stamm(2000)}.

\begin{theorem}\label{the: reduction of calf}
Let $\calf \subset \calg$ be families of subgroups of the group $\Gamma$. Let $\Lambda$ be a
ring with $\bbZ \subseteq \Lambda \subseteq \bbQ$ and  $N$ be an integer. Suppose for every $H\in \calg$ that
the assembly map induces for $n\le N$ an isomorphism
$$\Lambda \otimes_{\bbZ} \calh^{H}_n(\EGF{H}{H\cap \calf};\bfK(R?))  \to
\Lambda \otimes_{\bbZ}  \calh^{H}_n(H/H;\bfK(R?)),$$
where $H\cap \calf$ is the family of subgroups $K \subseteq H$ with $K \in \calf$.
Then the map
$$\Lambda \otimes_{\bbZ} H_n^{\Gamma}(\EGF{\Gamma}{\calf};\bfK(R?))  \to
\Lambda \otimes_{\bbZ}H_n^{\Gamma}(\EGF{\Gamma}{\calg};\bfK(R?))$$
is an isomorphism for $n \le N$. The analogous result is true for
$\bfL^{\langle - \infty \rangle}(R?)$ instead of $\bfK(R?)$.
\end{theorem}

In the sequel we will apply Theorem~\ref{the: reduction of calf} using the fact
that for an infinite cyclic group or an infinite dihedral group $H$ the map
$$\asmb \colon \calh^H_n(\underline{E}H;\bfK(R?)) ~ \to ~
\calh^K_n(H/H;\bfK(R?)) = K_n(RH)$$
is bijective for $n \in \bbZ$. This follows for the infinite cyclic group from the Bass--Heller-decomposition
and for the infinite dihedral group from Waldhausen
\cite[Corollary 11.5 and the following
Remark]{Waldhausen(1978b)} (see also \cite{Bartels-Lueck(2004)} and \cite[Section 2.2]{Lueck-Reich(2003b)}).

The Farrell--Jones Conjecture for algebraic $K$--theory for the trivial family
$\caltr$ consisting of the trivial subgroup only is true
for infinite cyclic groups and regular rings $R$ as coefficients. We conclude from
Theorem~\ref{the: reduction of calf} that for a regular ring $R$ the maps
\begin{multline*}
\calh^{p^{-1}(M_i)}_n(E(p^{-1}(M_i));\bfK(R?)) \\
\xrightarrow{\cong} ~
\calh^{p^{-1}(M_i)}_n(\EGF{p^{-1}(M_i)}{\calvcyc(K \cap p^{-1}(M_i))};\bfK(R?))
\end{multline*}
and
\begin{eqnarray*}\calh^G_n(EG;\bfK(R?)) & \xrightarrow{\cong} &
\calh^G_n(\EGF{G}{\calvcyc(K)};\bfK(R?))
\end{eqnarray*}
are bijective for all $n \in \bbZ$. Hence we obtain for a regular ring $R$
from the $G$--pushout \eqref{G-pushout for EGcalvcyc_f} an isomorphism
\begin{multline}
\bigoplus_{i \in I}
\calh^{p^{-1}(M_i)}_n\left(E(p^{-1}(M_i)) \to \EGF{p^{-1}(M_i)}{\calvcyc};\bfK(R?)\right)
\\
\xrightarrow{\cong}
\calh^G_n\left(EG \to \EGF{G}{\calvcyc_f};\bfK(R?)\right).
\label{calh^G_nleft(EG to EGF(G)(calvcyc_f);bfK(R?)right)}
\end{multline}
Let $\calvcyc_1$ be the family of virtually cyclic subgroups of $G$
whose intersection with $K$ is trivial. Since $\calvcyc$ is the union
$\calvcyc_f \cup\calvcyc_1$ and the intersection
$\calvcyc_f \cap \calvcyc_1$ is $\calfin$, we obtain a $G$--pushout
\begin{eqnarray}
&\comsquare{\underline{E}G}{}{\EGF{G}{\calvcyc_1}}{}{}
{\EGF{G}{\calvcyc_f}}{}{\EGF{G}{\calvcyc}} &
\label{G-pushout for calvcyc_f, calvcyc_1 and calvcyc}
\end{eqnarray}
The following conditions are equivalent for a virtually cyclic group $V$:
i.) $V$ admits an epimorphism to $\bbZ$ with finite kernel,
ii.) $H_1(V;\bbZ)$ is infinite, iii.) The center of $V$ is infinite.
A virtually cyclic subgroup does not satisfy these three equivalent conditions if and only if
it admits an epimorphism onto $D_{\infty}$ with finite kernel.

\begin{lemma} \label{lem: virtually cyclic subgroups of Q}
Any virtually cyclic subgroup of $Q$ is finite, infinite cyclic or isomorphic to $D_{\infty}$.
\end{lemma}
\begin{proof}
 Suppose that $V\subseteq Q$ is an infinite virtually cyclic subgroup.
Choose a finite normal subgroup $F \subseteq V$ such that $V/F$ is $\bbZ$ or $D_{\infty}$.
We have to show that $F$ is trivial. Suppose $F$ is not trivial. By assumption
there is a unique maximal finite subgroup $M \subseteq Q$ with $F \subseteq M$.
Consider $q \in N_GF$. Then $F \subseteq q^{-1}Mq \cap M$. This implies $q \in N_GM = M$.
Hence $N_GF$ is contained in the finite group $M$ what contradicts $V \subseteq N_GF$.
Hence $F$ must be trivial.
\end{proof}

Now we can prove Theorem~\ref{lem: result for algebraic $K$-theory}.
\begin{proof}
Lemma~\ref{lem: virtually cyclic subgroups of Q} implies that any
infinite subgroup appearing in $\calvcyc_1$ is an infinite cyclic
group or an infinite dihedral group. Hence Theorem~\ref{the:
reduction of calf} implies that $\calh^G_n\left(\underline{E}G \to
\EGF{G}{\calvcyc_1};\bfK(R?)\right)$ vanishes for $n \in \bbZ$. We
conclude from the $G$--pushout \eqref{G-pushout for calvcyc_f,
calvcyc_1 and calvcyc} that $\calh^G_n\left(\EGF{G}{\calvcyc_f}
\to \EGF{G}{\calvcyc};\bfK(R?)\right)$ vanishes for $n \in \bbZ$.
Now Theorem~\ref{lem: result for algebraic $K$-theory} follows
from \eqref{calh^G_nleft(EG to EGF(G)(calvcyc_f);bfK(R?)right)}.
\end{proof}

Now let us investigate what the results above imply for the middle and lower algebraic
$K$--theory with integral coefficients of the group $G = \Hei \rtimes \bbZ/4$ introduced
in Section~\ref{sec: The example Hei rtimes bbZ/4}
and $R = \bbZ$.
The Farrell--Jones Conjecture for algebraic K-theory is true
for $G$  and $R = \bbZ$ in the range
$n \le 1$ since $G$ is a discrete cocompact subgroup of the virtually connected Lie group
$\Hei(\bbR) \rtimes \bbZ/4$ (see \cite{Farrell-Jones(1993a)}). Each group $p^{-1}(M_i)$ is
virtually cyclic and satisfies the  Farrell--Jones Conjecture for algebraic K-theory
for trivial reasons. From Theorem~\ref{lem: result for algebraic $K$-theory}
we get for $n \le 1$ an isomorphism
$$\Wh_n(p^{-1}(M_0)) \bigoplus \Wh_n(p^{-1}(M_1)) \bigoplus
\Wh_n(p^{-1}(M_2)) \xrightarrow{\cong} \Wh_n(G),$$
which comes from the various inclusions of subgroups and the subgroups $M_0$, $M_1$ and $M_2$
of $Q$ have been introduced in Lemma~\ref{lem: maximal finite subgroups of Q}.
The Bass--Heller--Swan decomposition yields an isomorphism for any group $H$ 
\begin{eqnarray}
\hspace{-8mm} \Wh_n(H \times \bbZ) & \cong &
\Wh_{n-1}(H) \bigoplus \Wh_n(H) \bigoplus N\!K_n(\bbZ H) \bigoplus
\!NK_n(\bbZ H).
\label{Bass-Heller-Swan decomposition}
\end{eqnarray}
The groups $\Wh_n(\bbZ^k)$ and $\Wh_n(\bbZ/2 \times \bbZ^k)$ vanish for $n \le 1$
and $k \ge 0$. The groups $\Wh_n(\bbZ/4)$ are trivial for $n \le 1$.
References for these claims are given in the proof of \cite[Theorem 3.2]{Lueck-Stamm(2000)}.
The groups $\Wh_n(\bbZ/4 \times \bbZ^k)$ vanish for $n \le -1$ and $k \ge 0$.
This follows from \cite{Farrell-Jones(1995)}. Thus we get the following:

\begin{corollary} \label{cor: algebraic K-theory of G = Hei rtimes Z/4}
Let $G$ be the group $\Hei \rtimes \bbZ/4$ introduced
in Section~\ref{sec: The example Hei rtimes bbZ/4}. Then
\begin{eqnarray}
\Wh_n(G) & \cong & \left\{\begin{array}{lll}
N\!K_n(\bbZ [\bbZ/4]) \bigoplus N\!K_n(\bbZ [\bbZ/4])  & & \text { for  } n = 0,1;
\\
0 & & \text{ for } n \le -1.
\end{array}\right.
\label{explicite computation of Wh_n(G)}
\end{eqnarray}
where the isomorphism for $n = 0,1$ comes from the inclusions of
the subgroup $p^{-1}(M_0) = \langle t,z \rangle =\bbZ \times \bbZ/4$
into $G$ and the Bass--Heller--Swan
decomposition \eqref{Bass-Heller-Swan decomposition}.
\end{corollary}
Some information about $N\!K_n(\bbZ[\bbZ/4])$ is given
in \cite[Theorem 10.6 on page 695]{Bass(1968)}. Their
exponent divides $4^d$ for some natural number $d$.

\section{$L$--theory}
\label{sec: $L$-theory}

In this section we want to describe
what the methods above yield for the algebraic $L$--theory
provided that instead of the Baum--Connes
Conjecture the relevant version of the Farrell--Jones Conjecture for
algebraic $L$--theory (see \cite{Farrell-Jones(1993a)}) is true.

\begin{theorem}
\label{the: $L$-theory}
Let $G$ be the group appearing in \eqref{extension of groups} and assume
that conditions (M), (NM), and (T)  are satisfied.
Suppose that $G$ and all the groups $p^{-1}(M)$
for $M \subseteq Q$ maximal finite satisfy
the Farrell--Jones Conjecture for $L$--theory with coefficients in $R$. Then:
\begin{enumerate}

\item \label{the: $L$-theory: computation for underline E}
There is a long exact sequence of abelian groups
\begin{multline*}
\ldots \to
\calh_{n+1}(G\backslash\underline{E}G;\bfL^{\langle -\infty \rangle}(R))
\to \bigoplus_{i \in I} L_n^{\langle -\infty \rangle}(R[p^{-1}(M_i)])
\\
\to
\calh^G_n(\underline{E}G;\bfL^{\langle - \infty \rangle}(R?))  \bigoplus
\left(\bigoplus_{i \in I}
\calh_{n}(p^{-1}(M_i)\backslash\underline{E}p^{-1}(M_i);\bfL^{\langle -\infty \rangle}(R?))\right)
\\
 \to
\calh_n(G\backslash\underline{E}G;\bfL^{\langle -\infty \rangle}(R))\to
\bigoplus_{i \in I} L_{n-1}^{\langle -\infty \rangle}(R[p^{-1}(M_i)])
\to \ldots.
\end{multline*}
Let $\Lambda$ be a ring with $\bbZ \subseteq \Lambda \subseteq\bbQ$ such that the order
of each finite subgroup of $G$ is invertible in $\Lambda$. Then the long exact sequence above
reduces after applying $\Lambda \otimes_{\bbZ} -$ to short split exact sequences of
$\Lambda$--modules
\begin{multline*}
0 \to \bigoplus_{i \in I}
\Lambda \otimes_{\bbZ} L_n^{\langle -\infty \rangle}(R[p^{-1}(M_i)]) \to
\Lambda \otimes_{\bbZ} \calh^G_n(\underline{E}G;\bfL^{\langle - \infty \rangle}(R?))
\\  \bigoplus\left(\bigoplus_{i \in I}
 \Lambda \otimes_{\bbZ} \calh_{n}(p^{-1}(M_i)\backslash\underline{E}p^{-1}(M_i);
\bfL^{\langle -\infty \rangle}(R))\right)
\\
\to
\Lambda \otimes_{\bbZ}
\calh_{n}(G\backslash\underline{E}G;\bfL^{\langle -\infty \rangle}(R))  \to 0;
\end{multline*}

\item \label{the: $L$-theory: splitting due to Bartels}
Suppose for any virtually cyclic subgroup $V \subseteq G$ that
$K_{-i}(RV) = 0$ for sufficiently large $i$.
Then there is
a canonical isomorphism
\begin{multline*}
\calh^G_n(\underline{E}G;\bfL^{\langle -\infty \rangle}(\bbR?)) \bigoplus
\calh^G_n\left(\underline{E}G \to \EGF{G}{\calvcyc};\bfL^{\langle -\infty \rangle}(R?)\right)
\\
\xrightarrow{\cong}  L_n^{\langle - \infty  \rangle}(RG);
\end{multline*}

\item \label{the: $L$-theory: relative term}
We have
$$\calh^G_n\left(\underline{E}G \to \EGF{G}{\calvcyc};
\bfL^{\langle -\infty \rangle}(R?)\right)\left[\frac{1}{2}\right]  =  0.$$
\end{enumerate}
\end{theorem}
\begin{proof}
~\ref{the: $L$-theory: computation for underline E}\qua This is proved  in a
completely analogous way to
Theorem~\ref{the: modified Mayer-Vietoris sequence}.

\ref{the: $L$-theory: splitting due to Bartels}\qua  This follows from Theorem~\ref{the: Bartels results}.

\ref{the: $L$-theory: relative term}\qua  This follows from the next
Lemma~\ref{lem: general simplification of calvcyc}.
\end{proof}

\begin{lemma} \label{lem: general simplification of calvcyc}
Let $\Gamma$ be a group. Let $\calvcyc$ be the family of virtually cyclic subgroups of $\Gamma$
and $\calvcyc_{\bbZ}$ be the subfamily of $\calvcyc$ consisting of subgroups of $\Gamma$
which admit an epimorphism to $\bbZ$ with finite kernel. Let $\calf$ and $\calg$
be families of subgroups of $\Gamma$. If
$\calfin \subseteq \calf \subseteq \calg \subseteq \calvcyc_{\bbZ}$
holds, then
\begin{eqnarray*}
\calh^{\Gamma}_n\left(\EGF{\Gamma}{\calf} \to \EGF{\Gamma}{\calg};
\bfL^{\langle - \infty \rangle}(R?)\right) & = & 0.
\end{eqnarray*}
If
$\calfin \subseteq \calf \subseteq \calg\subseteq \calvcyc$
holds, then
\begin{eqnarray*}
\calh^{\Gamma}_n\left(\EGF{\Gamma}{\calf} \to \EGF{\Gamma}{\calg};
\bfL^{\langle - \infty \rangle}(R?)\right)\left[\frac{1}{2}\right] & = & 0.
\end{eqnarray*}
\end{lemma}
\begin{proof}
We conclude from Theorem~\ref{the: reduction of calf} that it suffices
to show for a virtually cyclic group $V$, which admits an epimorphism to
$\bbZ$, that the map
\begin{eqnarray*}
\calh^{V}_n(\underline{E}V;
\bfL^{\langle - \infty \rangle}(R?)) & \to & L_n^{\langle - \infty \rangle}(RV)
\end{eqnarray*}
is bijective and for a virtually cyclic group $V$, which admits an epimorphism
to $D_{\infty}$, that the map above is bijective after inverting two.

We begin with the case where $V = F \rtimes_{\phi} \bbZ$ for an automorphism
$\phi \colon F \to F$ of a finite group $F$.
There is a long exact sequence
which can be derived from \cite{Ranicki(1973b)} and \cite{Ranicki(1973c)}:
\begin{multline*}
\ldots \to L_n^{\langle - \infty \rangle}(RF)
\xrightarrow{\id -  L_n^{\langle - \infty \rangle}(R\phi)}
L_n^{\langle - \infty \rangle}(RF) \to L_n^{\langle - \infty \rangle}(RV) \to
\\
\to L_{n-1}^{\langle - \infty \rangle}(RF)
\xrightarrow{\id -  L_{n-1}^{\langle - \infty \rangle}(R\phi)}
L_{n-1}^{\langle - \infty \rangle}(RF) \to \ldots
\end{multline*}
Since $\bbR$ with the action of $V$ coming from the epimorphism to $\bbZ$ and the
action of $\bbZ$ by translation is a model for $\underline{E}V$, we also obtain a long exact
Mayer--Vietoris sequence:
\begin{multline*}
\ldots \to L_n^{\langle - \infty \rangle}(RF)
\xrightarrow{\id -  L_n^{\langle - \infty \rangle}(R\phi)}
L_n^{\langle - \infty \rangle}(RF) \to
\calh^V_n(\underline{E}V;\bfL^{\langle - \infty \rangle}(R?)) \to
\\
\to L_{n-1}^{\langle - \infty \rangle}(RF)
\xrightarrow{\id -  L_{n-1}^{\langle - \infty \rangle}(R\phi)}
L_{n-1}^{\langle - \infty \rangle}(RF) \to \ldots
\end{multline*}
These two sequences are compatible with the assembly map
$$\asmb \colon \calh^V_n(\underline{E}V;\bfL^{\langle - \infty \rangle}(R?)) \to
\calh^V_n(V/V;\bfL^{\langle - \infty \rangle}(R?)) = L_n^{\langle -\infty\rangle}(RV),
$$
which must be an isomorphism by the Five-Lemma.

Suppose that $V$ admits an epimorphism onto $D_{\infty} = \bbZ/2 \ast \bbZ/2$.
Then we can write $V$ as an amalgamated product $F_1 \ast_{F_0} F_2$ for
finite groups $F_1$ and $F_2$ and a common subgroup $F_0$.
We can think of $V$ as a graph of groups associated to a
segment and obtain an action without inversions on a tree
which yields a $1$--dimensional model for $\underline{E}V$ with two equivariant
$0$--cells $V/F_1$ and $V/F_2$ and one equivariant one-cell $V/F_0 \times D^1$
(see \cite[\S 5]{Serre(1980)}). The associated long Mayer--Vietoris sequence looks like:
\begin{multline*}
\ldots \to L^{\langle - \infty \rangle}_n(RF_0)
\to L^{\langle - \infty \rangle}_n(RF_1) \bigoplus L^{\langle - \infty \rangle}_n(RF_2)
\to
\calh^V_n(\underline{E}V;\bfL^{\langle - \infty \rangle}(R?))
\\
\to L^{\langle - \infty \rangle}_{n-1}(RF_0)
\to L^{\langle - \infty \rangle}_{n-1}(RF_1) \bigoplus L^{\langle - \infty \rangle}_{n-1}(RF_2)
\to \ldots
\end{multline*}
There is a corresponding exact sequence, where
$\calh^V_n(\underline{E}V;\bfL^{\langle - \infty \rangle}(R?))$ is replaced by
$L^{\langle  -\infty \rangle}_n(RV)$ and additional $\UNil$--terms occur
which vanish after inverting two (see for $\bbZ \subseteq R \subseteq \bbQ$ \cite[Corollary 6]{Cappell(1974c)}
or see \cite[Remark 8.7]{Ranicki(1995b)} and \cite{Ranicki(1973b)}).
Now a Five-Lemma argument proves the claim.
\end{proof}

\begin{theorem}
\label{the: $L$-theory for Q without elements of order two}
Let $G$ be the group appearing in \eqref{extension of groups} and assume
that conditions (M), (NM), and (T)  are satisfied. Suppose that $Q$ contains no element
of order $2$. Suppose that $G$ and all the groups $p^{-1}(M)$
for $M \subseteq Q$ maximal finite satisfy
the Farrell--Jones Conjecture for $L$--theory with coefficients in $R$.
Then there is a long exact sequence of abelian groups:
\begin{multline*}
\ldots \to
\calh_{n+1}(G\backslash\underline{E}G;\bfL^{\langle -\infty \rangle}(R?)) \to 
\bigoplus_{i \in I} L_n^{\langle -\infty \rangle}(R[p^{-1}(M_i)])
\\
\to
L^{\langle - \infty \rangle}_n(RG)  \bigoplus
\left(\bigoplus_{i \in I}
\calh_{n}(p^{-1}(M_i)\backslash\underline{E}p^{-1}(M_i);\bfL^{\langle -\infty \rangle}(R?))\right)
\\
\to \calh_{n}(G\backslash\underline{E}G;\bfL^{\langle -\infty \rangle}(R?)) 
\to \ldots
\end{multline*}
Let $\Lambda$ be a ring with $\bbZ \subseteq \Lambda \subseteq\bbQ$ such that the order
of each finite subgroup of $G$ is invertible in $\Lambda$. Then the long exact sequence above
reduces after applying $\Lambda \otimes_{\bbZ} -$ to short split exact sequences of
$\Lambda$--modules
\begin{multline*}
0 \to \bigoplus_{i \in I} \Lambda \otimes_{\bbZ} L_n^{\langle -\infty \rangle}(R[p^{-1}(M_i)]) \to
\\
\Lambda \otimes_{\bbZ} L^{\langle - \infty \rangle}_n(RG)
\bigoplus  \left(\bigoplus_{i \in I}
\Lambda \otimes_{\bbZ} \calh_{n}(p^{-1}(M_i)\backslash\underline{E}p^{-1}(M_i);
\bfL^{\langle -\infty \rangle}(R))\right)
\\
\to
\Lambda \otimes_{\bbZ}
\calh_{n}(G\backslash\underline{E}G;\bfL^{\langle -\infty \rangle}(R))  \to 0.
\end{multline*}
\end{theorem}
\begin{proof}
Because of Theorem~\ref{the: $L$-theory} and
Lemma~\ref{lem: general simplification of calvcyc}
it suffices to prove 
that $V$ admits an epimorphism
to $\bbZ$ for an infinite virtually cyclic subgroup $V \subset G$. 
If $V \cap K$ is trivial, then $V$ is an infinite virtually cyclic subgroup of
$Q$ and hence isomorphic to $\bbZ$ by Lemma~\ref{lem: virtually cyclic subgroups of Q}.
Suppose that $V \cap K$ is non-trivial. Then $V$ can be written as an extension
$1 \to K \cap V \to V \to p(V) \to 1$ for a finite subgroup $p(V) \subseteq Q$.
The group $K \cap V$ is infinite cyclic and $p(V)$ must have odd order.
Hence $V$ contains a central infinite cyclic subgroup. This implies that
$V$ admits an epimorphism to $\bbZ$.
\end{proof}

From now on we assume that $Q$ is an extension
$1 \to \bbZ^n \to Q \to F \to 1$ for a finite group $F$ such that
the conjugation action of $F$ on $\bbZ^n$ is free outside the origin.

Let $V \subseteq Q$ be infinite virtually cyclic.  Either $F$ is of odd order or
$F$ has a unique element $f_2$ of order two \cite[Lemma 6.2]{Lueck-Stamm(2000)}.
Because of Lemma~\ref{lem: virtually cyclic subgroups of Q} either $V$ is infinite cyclic or
$V$ is isomorphic to $D_{\infty}$ and each element $v_2 \in V$ of order two is mapped
under $Q \to F$ to the unique element $f_2 \in F$ of order two in $F$.

Suppose that $V \subseteq Q$ is isomorphic to $D_{\infty}$.
Then $V$ contains at least one element
$v_2 \in V$ of order two. Any other element of order two is of the shape
$v_2u$ for $u \in V \cap \bbZ^n$. Hence $V$ is the subgroup
$\langle v_2, V \cap \bbZ^n \rangle$ generated by $v_2$ and the infinite
cyclic group $V \cap \bbZ^n$ regardless which element $v_2 \in V$ of order two we choose.
For any infinite cyclic subgroup $C \subseteq \bbZ^n$ let $C_{\max}$ be the kernel of the projection
$\bbZ^n \to (\bbZ^n/C)/\tors(\bbZ^n/C)$. This is the maximal infinite cyclic subgroup
of $\bbZ^n$ which contains $C$. Define $V_{\max} \subseteq Q$ to be
$$V_{\max} := \langle v_2, (V \cap \bbZ^n)_{\max}\rangle.$$
The subgroup $V_{\max}$  is
isomorphic to $D_{\infty}$ and satisfies $V \subseteq V_{\max}$. Moreover, it is a maximal virtually
cyclic subgroup, ie, $V_{\max} \subseteq W$ for a virtually cyclic subgroup
$W \subseteq Q$ implies $V_{\max} = W$. Let $V \subseteq W$ be virtually cyclic subgroups
of $Q$ such that $V \cong D_{\infty}$. Then $W \cong D_{\infty}$ and $V_{\max} =
W_{\max}$. Hence each virtually cyclic subgroup $V$ with $V \cong D_{\infty}$ is contained
in a unique maximal virtually cyclic subgroup of $Q$ isomorphic to $D_{\infty}$, namely $V_{\max}$.

Next we show $N_GV = V$ if $V$ is a subgroup of $Q$ with
$V \cong D_{\infty}$ and $V = V_{\max}$. Let $v_2 \in V$ be an element of order two.
Consider an element $q \in N_GV$. Then the conjugation action of $q$ on $\bbZ^n$
sends $V \cap \bbZ^n$ to itself. Hence the conjugation action of $q^2$ on $\bbZ^n$ induces
the identity on $V \cap \bbZ^n$. This implies that $q$ is mapped under $Q \to F$
to the unit element or the unique element of order two $f_2$.
Hence $q$ is of the shape $u$ or $v_2u$ for some $u \in \bbZ^n$.
Since $(v_2u)v_2(v_2u)^{-1} = v_2uv_2u^{-1}v_2 = v_2u^2$ and
$uv_2u^{-1} = v_2u^2$, we conclude $v_2u^2 \in V$. This implies that
$u^2 \in V\cap \bbZ^n$. Since $V = V_{\max}$, we get $u \in  V$ and hence
$q\in V$.

Let $J$ be a complete system of representatives $V$ for the
set of conjugacy classes $(V)$ of subgroups $V \subseteq Q$ with
$V \cong D_{\infty}$ and $V = V_{\max}$. In the sequel let
$\calicof$ be the set of subgroups $H$ with are infinite cyclic or finite.
By attaching equivariant cells we construct a model for $\EGF{Q}{\calicof}$ which contains
$\coprod_{V \in J} Q \times_V \EGF{V}{\calicof}$
as $Q$--$CW$--subcomplex. Define a $Q$--$CW$--complex $\EGF{Q}{\calvcyc}$ by the
$Q$--pushout
\begin{eqnarray}
& \comsquare{\coprod_{V \in J} Q \times_V \EGF{V}{\calicof}}
{}{\EGF{Q}{\calicof}}
{u_5}{f_5}
{\coprod_{V \in J} Q/V}
{}{\EGF{Q}{\calvcyc}}
& \label{Q-pushout for EGF Q calvcyc}
\end{eqnarray}
where the map $u_5$ is the obvious projection and the upper horizontal arrow is
the inclusion.

We have to show that $\EGF{Q}{\calvcyc}$ is a model for the classifying space for the
family $\calvcyc$ of virtually cyclic subgroups of $Q$.
Obviously all its isotropy groups belong to $\calvcyc$.
Let $H \subseteq Q$ be a virtually cyclic group with $H \in \calicof$.
Choose a map of sets $s \colon Q/V \to Q$ such that its composition with the projection
$Q \to Q/V$ is the identity. For any $V$--space $X$, there is a homeomorphism
$$\left(Q \times_V X\right)^H \xrightarrow{\cong}
\coprod_{\substack{qV \in Q/V\\s(qV)^{-1}Hs(qV) \subseteq V}} X^{s(qV)^{-1}Hs(qV)},
\hspace{5mm} (q,x)  \mapsto s(qV)^{-1}qx,$$
whose inverse sends $x$ of the summand $X^{s(qV)^{-1}Hs(qV)}$
belonging to $qV \in Q/V$ with $s(qV)^{-1}Hs(qV) \subseteq V$ to
$(s(qV),x)$. Since $\EGF{Q}{\calicof}^{s(qV)^{-1}Hs(qV)}$
and $(V/V)^{s(qV)^{-1}Hs(qV)}$ are contractible for each $qV \in Q/V$ with
$s(qV)^{-1}Hs(qV) \subseteq V$, the map $u_5^H$ is a homotopy equivalence.
Hence $f_5^H$ is a homotopy equivalence. The space $\EGF{Q}{\calicof}^H$ is contractible.
Therefore $\EGF{Q}{\calvcyc}^H$ is contractible.
Let $H \subseteq Q$ be a virtually cyclic group with $H \not\in \calicof$. Then
\begin{eqnarray*}
\left(\coprod_{V \in J} Q/V\right)^H &  = & \pt;
\\
\left(Q \times_V \EGF{V}{\calicof}\right)^H
& = &
\emptyset;
\\
\EGF{Q}{\calicof}^H
& = &
\emptyset.
\end{eqnarray*}
This implies that $\EGF{Q}{\calvcyc}^H = \pt$ is contractible.

Recall that $\calvcyc_f$ is the family of virtually cyclic
subgroups of $G$ whose image under $p\colon G \to Q$ is finite
and $\calvcyc_1$ is the family of virtually cyclic
subgroups of $G$ whose intersection with $K = \ker(p)$ is trivial.
Let $\calvcyc_{icof}$ be the family of subgroups of $G$
whose image under $p \colon G \to Q$ is contained in $\calicof$.
If we cross the $Q$--pushout \eqref{Q-pushout for EGF Q calvcyc} with
$\EGF{G}{\calvcyc}$, we obtain the $G$--pushout:
\begin{eqnarray*}
& \comsquare{\coprod_{V \in J} G \times_{p^{-1}(V)} \EGF{p^{-1}(V)}{\calicof}}
{}{\EGF{G}{\calvcyc_{icof}}}
{u_6}{f_6}
{\coprod_{V \in J} G \times_{p^{-1}(V)} \EGF{p^{-1}(V)}{\calvcyc}}
{}{\EGF{G}{\calvcyc}}
&
\end{eqnarray*}
Because of Lemma~\ref{lem: general simplification of calvcyc}
this $G$--pushout induces isomorphisms for $n \in \bbZ$ 
\begin{multline*}
\bigoplus_{V \in J} \calh^{p^{-1}(V)}_n\left(\underline{E}p^{-1}(V) \to \EGF{p^{-1}(V)}{\calvcyc};
\bfL^{\langle -\infty \rangle}(R?)\right)
\\
\xrightarrow{\cong}
\calh^G_n\left(\EGF{G}{\calvcyc_f} \to \EGF{G}{\calvcyc}; \bfL^{\langle -\infty \rangle}(R?)\right).
\end{multline*}
We now conclude from Lemma~\ref{lem: general simplification of calvcyc}:

\begin{lemma} \label{lem: Z^n to Q to F and condition and calvcyc_f}
Let $1 \to \bbZ^n \to Q \to F \to 1$ be an extension such that
the conjugation action of $F$ on $\bbZ^n$ is free outside the origin.
Let $ 1\to K \xrightarrow{i}  G \xrightarrow{p} Q \to 1$ be an extension.
Suppose that for any virtually cyclic group
$V \subseteq G$ with $p(V)$ finite there exists an epimorphism $V \to \bbZ$.
(This condition is satisfied if $K$ is abelian and contained  in the center of $G$.)
Then there is an isomorphism
\begin{multline*}
\bigoplus_{V \in J} \calh^{p^{-1}(V)}_n\left(\underline{E}p^{-1}(V) \to \EGF{p^{-1}(V)}{\calvcyc};
\bfL^{\langle -\infty \rangle}(R?)\right)
\\
\xrightarrow{\cong}
\calh^G_n\left(\underline{E}G \to \EGF{G}{\calvcyc}; \bfL^{\langle -\infty \rangle}(R?)\right).
\end{multline*}
\end{lemma}

Next we apply Theorem~\ref{the: $L$-theory} and Lemma~\ref{lem:
Z^n to Q to F and condition and calvcyc_f} to the special example
$G = \Hei \rtimes \bbZ/4$ introduced in Section~\ref{sec: The
example Hei rtimes bbZ/4}. We begin with constructing an explicit
choice for $J$ and determining the preimages $p^{-1}(V)$ for $V
\in J$. Recall that $J$ is a complete system of representatives of
the conjugacy classes $(V)$ of subgroups $V \subseteq Q$ with $V
\cong D_{\infty}$ and $V = V_{\max}$. Let $IC(\bbZ^2)$ be the set
of infinite cyclic subgroups $L$ of $\bbZ^2$. Any subgroup $V
\subseteq Q$ with $V \cong D_{\infty}$ can be written as $V =
\langle v_2,V \cap \bbZ^2\rangle$ for $v_2 \in V$ any element of
order two. Hence we can write $V = \langle t^2a ,L\rangle$ for $L
\in IC(\bbZ^2)$ and $a \in \bbZ^2$. We have $\langle t^2 a,L
\rangle = \langle t^2 a',L' \rangle$ if and only if $L = L'$ and
$a -a' \in L = L'$. We have $V = V_{\max}$ for $V = \langle t^2a
,L\rangle$ if and only if $L \subset \bbZ^2$ is maximal.

Let $IC^+(\bbZ)$ be the subset for which $L \in IC(\bbZ^2)$ meets
$\{(n_1,n_2) \in \bbZ^2 \mid n_1 \ge 0, n_2 > 0\}$. The $\bbZ/4$--action on $\bbZ^2$
induces a $\bbZ/2$--action on $L_1(\bbZ^2)$ by sending $L$ to $i\cdot L$.
Notice that $IC^+(\bbZ^2)$ is a fundamental domain for this action, ie,
$IC(\bbZ^2)$ is the  disjoint union of $IC^+(\bbZ^2)$ and its image under this
involution. We claim that a complete system of representatives of conjugacy classes (V)
of subgroups $V$ of $Q = \bbZ^2 \rtimes \bbZ/4$ with $V \cong
D_{\infty}$ is
\begin{eqnarray}
& \begin{array}{lll}
\langle t^2,(n_1,n_2)\rangle & & n_1 \text{ even};
\\
\langle t^2(0,1),(n_1,n_2)\rangle & & n_1 \text{ even};
\\
\langle t^2,(n_1,n_2)\rangle & & n_2 \text{ even};
\\
\langle t^2(1,0),(n_1,n_2)\rangle & & n_2 \text{ even};
\\
\langle t^2,(n_1,n_2)\rangle & & n_1 \text { and } n_2 \text{ odd};
\\
\langle t^2(1,0),(n_1,n_2) \rangle & & n_1 \text { and } n_2 \text{ odd},
\end{array} &
\label{list of the V-s}
\end{eqnarray}
where $(n_1,n_2)$ runs through $IC^+ = \{(n_1,n_2) \in \bbZ^2 \mid n_1 > 0,n_2 \ge 0,
(n_1,n_2) = 1\}$. This follows from the computations
\begin{eqnarray*}
(m_1,m_2)^{-1}(t^2(n_1,n_2))(m_1,m_2) & = & t^2(n_1 + 2m_1,n_2 + 2m_2);
\\
t(t^2(n_1,n_2))t^{-1} & = & t^2(-n_2,n_1);
\\
t^2(t^2(n_1,n_2))(t^2)^{-1} & = & t^2(-n_1,-n_2).
\end{eqnarray*}
Now we list the preimages $p^{-1}(V)$ of these
subgroups above and determine their isomorphism type. We claim that
they can be described by the following generators
\begin{eqnarray}
\hspace{-3mm} \begin{array}{lll}
\langle t^2,t^2(n_1,n_1n_2/2,n_2),(0,1,0)\rangle & & n_1 \text{ even};
\\
\langle t^2(0,0,1),t^2(n_1,n_1n_2/2,n_2),(0,1,0)\rangle & & n_1 \text{ even};
\\
\langle t^2,t^2(n_1,n_1n_2/2,n_2),(0,1,0)\rangle & & n_2 \text{ even};
\\
\langle t^2(1,0,0),t^2(n_1,n_1n_2/2,n_2),(0,1,0)\rangle & & n_2 \text{ even};
\\
\langle t^2,t^2(2n_1,2n_1n_2,2n_2),
(n_1, \frac{n_1n_2+1}{2},n_2)\rangle & & n_1 \text { and } n_2 \text{ odd};
\\
\langle t^2(1,0,0),t^2(2n_1+1,2n_1n_2+n_2,2n_2),(n_1,
\frac{n_1n_2+1}{2},n_2)
\rangle & & n_1 \text { and } n_2 \text{ odd},
\end{array}
& &
\label{list of the p^(-1)(V)-s}
\end{eqnarray}
where $(n_1,n_2)$ runs through $IC^+ = \{(n_1,n_2) \in \bbZ^2 \mid n_1 > 0,n_2 \ge 0,
(n_1,n_2) = 1\}$. This is obvious for the first four groups and follows for the
last two groups from the computation
\begin{eqnarray*}
(n_1,\frac{n_1n_2 +1}{2},n_2)^2 & = & (2n_1,2n_1n_2,2n_2) \cdot (0,1,0);
\\
(2n_1+1,2n_1n_2 + n_2,2n_2) & = & (1,0,0) \cdot (2n_1,2n_1n_2,2n_2).
\end{eqnarray*}
The first four groups are isomorphic to $D_{\infty}
\times \bbZ$ and the last two are isomorphic to the semi-direct product
$D_{\infty} \rtimes_a \bbZ$ with respect to the automorphism $a$ of $D_{\infty} = \bbZ/2\ast
\bbZ/2 = \langle s_1,s_2 \mid s_2^2 = s_2^2 = 1\rangle$ which send
$s_1$ to $s_2$ and $s_2$ to $s_1$. For the first four groups there are explicit isomorphisms
from $D_{\infty} \times \bbZ = \langle s_1,s_2,z \mid s_1^2 = s_2^2 =
[s_1,z] = [s_2,z] = 1 \rangle$ which send $s_1, s_2, z$ to the three
generators appearing in the presentation above. Similarly for the last two groups there are
explicit isomorphisms from $D_{\infty}\rtimes_a \bbZ = \langle
s_1,s_2,z \mid s_1^2 = s_2^2 = 1, z^{-1}s_1z = s_2\rangle$ which send
$s_1, s_2, z$ to the three
generators appearing in the presentation above. We leave it to the reader to check that
these generators appearing in the presentation above do satisfy the required relations.

Next we compute the groups
$\calh_n^{D_{\infty} \times \bbZ}\left(\underline{E}D_{\infty} \times \bbZ \to
\EGF{D_{\infty} \times \bbZ}{\calvcyc};\bfL^{-\infty}(\bbZ?)\right)$ and
$\calh_n^{D_{\infty} \rtimes_a \bbZ}\left(\underline{E}D_{\infty} \rtimes_a \bbZ \to
\EGF{D_{\infty} \rtimes \bbZ}{\calvcyc};\bfL^{-\infty}(\bbZ?)\right)$.
There is an obvious model for $\underline{E}D_{\infty}$,
namely $\bbR$ with the trivial $\bbZ$--action and the action of $D_{\infty} = \bbZ \rtimes_a
\bbZ/2$, which comes from the $\bbZ$--action by translation and the $\bbZ/2$--action given
by $-\id_{\bbR}$. From this we obtain an exact sequence
\begin{multline*}
0 \to L_n^{\langle - \infty\rangle }(R) \xrightarrow{i}
L_2^{\langle -\infty \rangle}(R[\bbZ/2]) \bigoplus
L_2^{\langle -\infty\rangle }(R[\bbZ/2]) 
\\
\xrightarrow{f}
\calh^{D_{\infty}}_n(\underline{E}D_{\infty};\bfL^{-\infty}(R?)) \to 0
\end{multline*}
such that the composition of $f$ with the obvious map
$$\calh^{D_{\infty}}_n(\underline{E}D_{\infty};\bfL(R?)) \to
\calh^{D_{\infty}}_n(\pt;\bfL(R?) = L^{\langle -\infty \rangle}_n(R[D_{\infty}])$$
is given by the two obvious inclusions
$\bbZ/2 \to D_{\infty} = \langle s_1,s_2 \mid s_1 = s_2 ^2 = 1\rangle$.
Thus we obtain an isomorphism
\begin{eqnarray}
\hspace{-9mm} \calh^{D_{\infty}}_n\left(\underline{E}D_{\infty} \to
\EGF{D_{\infty}}{\calvcyc};
\bfL^{\langle -\infty \rangle}(R?)\right)
& = & \UNil_n(\bbZ/2 \ast \bbZ/2;R),
\label{relative term for D_(infty)}
\end{eqnarray}
where $\UNil_n(\bbZ/2 \ast \bbZ/2;R)$ is the $\UNil$--term appearing in the
short split exact sequence
\begin{multline*}
0 \to L^{\langle - \infty\rangle}_n(R) \to L^{\langle - \infty\rangle}_n(R[\bbZ/2]) \oplus
L_n^{\langle - \infty\rangle}(R[\bbZ/2]) \oplus  \UNil_n(\bbZ/2 \ast \bbZ/2;R)
\\
\to L^{\langle - \infty\rangle}_n(R[\bbZ/2 \ast \bbZ/2]) \to 0
\end{multline*}
due to Cappell~\cite[Theorem 10]{Cappell(1974c)}.
For the computation of these terms $\UNil_n(\bbZ/2 \ast \bbZ/2;R)$
we refer to \cite{Banagl-Ranicki(2003)},
\cite{Connolly-Davis(2004)} and \cite{Connolly-Ranicki(2003)}. They have
exponent four and they are either trivial or are infinitely generated as abelian groups.

We can take as model for $\underline{E}(D_{\infty} \times \bbZ)$ the product
$\underline{E}D_{\infty} \times \bbR$, where $\bbZ$ acts on $\bbR$ by translation.
We get from  \eqref{relative term for D_(infty)} and Lemma~\ref{lem: general simplification of calvcyc}
using a Mayer--Vietoris argument an isomorphism
\begin{multline}
\calh^{\bbZ \times D_{\infty}}_n\left(\underline{E}(D_{\infty} \times \bbZ)\to
\EGF{D_{\infty} \times \bbZ}{\calvcyc};
\bfL^{\langle -\infty \rangle}(R?)\right)
\\
\cong
\calh^{\bbZ \times D_{\infty}}_n\left(\underline{E}D_{\infty} \times \bbR\to
\EGF{D_{\infty}}{\calvcyc}\times \bbR;
\bfL^{\langle -\infty \rangle}(R?)\right)
\\
\cong
\calh^{D_{\infty}}_n\left(\underline{E}D_{\infty} \times S^1\to
\EGF{D_{\infty}}{\calvcyc}\times S^1;
\bfL^{\langle -\infty \rangle}(R?)\right)
\\
\cong  \UNil_n(\bbZ/2 \ast \bbZ/2;R) \bigoplus \UNil_{n-1}(\bbZ/2 \ast \bbZ/2;R).
\label{relative term for D_(infty) times bbZ}
\end{multline}
Next we investigate $D_{\infty} \rtimes_a \bbZ$. Let
$\calvcyc(D_{\infty})$  be the family of virtually cyclic subgroups of
$D_{\infty}\rtimes_a \bbZ$ which lie in $D_{\infty}$ and let
$\calvcyc_f$ be the family of virtually cyclic subgroups of
$\calvcyc(D_{\infty})$ whose intersection with $D_{\infty}$ is
finite. Then the family $\calvcyc$ of virtually cyclic subgroups of
$D_{\infty}\rtimes_a \bbZ$ is the union of $\calvcyc(D_{\infty})$ and
$\calvcyc_f$ and the family $\calfin $ of finite subgroups of
$D_{\infty}\rtimes_a \bbZ$ is the intersection of $\calvcyc(D_{\infty})$ and
$\calvcyc_f$. Hence we get a pushout of $D_{\infty}\rtimes_a
\bbZ$--spaces
$$\comsquare{\underline{E}(D_{\infty}\rtimes_a \bbZ)}{}
{\EGF{D_{\infty}\rtimes_a \bbZ}{\calvcyc(D_{\infty})}}
{}{}
{\EGF{D_{\infty}\rtimes_a \bbZ}{\calvcyc_f}}{}
{\EGF{D_{\infty}\rtimes_a \bbZ}{\calvcyc}}
 $$
Any finite subgroup of $D_{\infty}$ is trivial or isomorphic to
$\bbZ/2$. Any group $V$ which can be written as extension $1 \to
\bbZ/2 \to V \to \bbZ \to 1$ is isomorphic to $\bbZ/2 \times
\bbZ$. Hence any infinite group $V$ occurring in $\calvcyc_f$ is
isomorphic to $\bbZ$ or $\bbZ/2 \times \bbZ/2$. We conclude from
Lemma~\ref{lem: general simplification of calvcyc} and the
$D_{\infty}\rtimes_a \bbZ$--pushout above
\begin{multline*}
\calh^{D_{\infty} \rtimes_a \bbZ)}_n\left(
\EGF{D_{\infty}\rtimes_a \bbZ}{\calvcyc(D_{\infty})} \to
\EGF{D_{\infty}\rtimes_a \bbZ}{\calvcyc};\bfL^{\langle - \infty \rangle}(R?)\right)
\\  \cong
\calh^{D_{\infty} \rtimes_a \bbZ)}_n\left(
\underline{E}(D_{\infty}\rtimes_a \bbZ) \to \EGF{D_{\infty}\rtimes_a \bbZ}{\calvcyc_f};
\bfL^{\langle - \infty \rangle}(R?)\right)   \cong  0.
\end{multline*}
Hence we get an isomorphism
\begin{multline*}
\calh^{D_{\infty} \rtimes_a \bbZ)}_n\left(
\underline{E}(D_{\infty}\rtimes_a \bbZ) \to \EGF{D_{\infty}\rtimes_a \bbZ}{\calvcyc};
\bfL^{\langle - \infty \rangle}(R?)\right)
\\  \cong
\calh^{D_{\infty} \rtimes_a \bbZ)}_n\left(
\underline{E}(D_{\infty}\rtimes_a \bbZ) \to \EGF{D_{\infty}\rtimes_a \bbZ}{\calvcyc(D_{\infty})};
\bfL^{\langle - \infty \rangle}(R?)\right).
\end{multline*}
We can take as model for $\underline{E}(D_{\infty} \rtimes_a \bbZ)$ the to both sides
infinite mapping telescope of the $\left(a\colon D_{\infty} \to D_{\infty}\right)$--equivariant map
$\underline{E}a \colon \underline{E}D_{\infty} \to \underline{E}D_{\infty}$ with the
$D_{\infty}\rtimes_a \bbZ$--action for which $\bbZ$ acts by shifting the telescope to the right.
A model for $\EGF{D_{\infty}\rtimes_a \bbZ}{\calvcyc(D_{\infty})}$ is the
to both sides
infinite mapping telescope of the $\left(a\colon D_{\infty} \to D_{\infty}\right)$--equivariant map
$\pt \to \pt$. Of course this is the same as $\bbR$ with the $D_{\infty}
\rtimes_a \bbZ$--action, for which $D_{\infty}$ acts trivially and $\bbZ$ by translation.
The long Mayer--Vietoris sequence together with \eqref{relative term for D_(infty)}
yields a long exact sequence:
\begin{multline}
\ldots \to \UNil_n(\bbZ/2 \ast \bbZ/2;R) \xrightarrow{\id - \UNil_n(a)} \UNil_n(\bbZ/2 \ast
\bbZ/2;R)
\\
\to \calh^{D_{\infty} \rtimes_a \bbZ}_n\left(
\underline{E}(D_{\infty}\rtimes_a \bbZ) \to \EGF{D_{\infty}\rtimes_a \bbZ}{\calvcyc};
\bfL^{\langle - \infty \rangle}(R?)\right)
\\
\to \UNil_{n-1}(\bbZ/2 \ast \bbZ/2;R) \xrightarrow{\id - \UNil_{n-1}(a)} \UNil_{n-1}(\bbZ/2 \ast
\bbZ/2;R) \to \ldots
\label{relative term for D rtimes Z}
\end{multline}
The homomorphism $\UNil_{n-1}(a)$ has been analyzed in \cite{Brookman-Davis-Khan(2005)}.

\begin{theorem} \label{the: L-theory of hei rtimes Z/4}
Let $G$ be the group $\Hei \rtimes \bbZ/4$
introduced in Section~\ref{sec: The example Hei rtimes bbZ/4}. Then

\begin{enumerate}

\item  \label{the: L-theory of hei rtimes Z/4: underline(E)}
There is a short exact sequence which splits after inverting $2$
\begin{multline*}
0 \to L_n^{\langle -\infty \rangle}(\bbZ) \bigoplus
\widetilde{L}_n^{\langle -\infty \rangle}(\bbZ[\bbZ/2]) \bigoplus
\widetilde{L}^{\langle -\infty \rangle}_{n-1}(\bbZ[\bbZ/2])
\\
\bigoplus \widetilde{L}^{\langle -\infty \rangle}_n(\bbZ[\bbZ/4]) \bigoplus
\widetilde{L}^{\langle -\infty \rangle}_{n-1}(\bbZ[\bbZ/4])
\\
\xrightarrow{j} \calh^G(\underline{E}G;\bfL^{\langle -\infty \rangle}(\bbZ))
\to L_{n-3}^{\langle -\infty \rangle}(\bbZ) \to 0;
\end{multline*}

\item  \label{the: L-theory of hei rtimes Z/4: splitting for fin and relative term}
There is for $n \in \bbZ$ an isomorphism
$$
\calh^G_n(\underline{E}G;\bfL^{\langle -\infty \rangle}(\bbZ))
\bigoplus
\calh^G_n\left(\underline{E}G \to \EGF{G}{\calvcyc};\bfL^{\langle -\infty \rangle}(\bbZ?)\right)
\xrightarrow{\cong}  L_n^{\langle - \infty \rangle}(\bbZ G);
$$

\item \label{the: L-theory of hei rtimes Z/4: relative term}
Let $IC^+$ be the set $\{(n_1,n_2) \in \bbZ^2 \mid n_1 > 0, n_2 \ge 0, (n_1,n_2) = 1\}$.
Then there is an isomorphism
\begin{multline*} \hspace{-4mm}\Bigg(
\bigoplus_{\substack{(n_1,n_2) \in IC^+,\\ n_1 \text{ or } n_2 \text{ even}}}
\bigoplus_{i=1}^4 \left(\UNil_n(\bbZ/2 \ast \bbZ/2;\bbZ)
\bigoplus \UNil_{n-1}(\bbZ/2 \ast \bbZ/2;\bbZ)\right)\Bigg) \bigoplus
\\\hspace{-4mm}
\Bigg(
\bigoplus_{\substack{(n_1,n_2) \in IC^+,\\ n_1 \text{ and } n_2 \text{ odd}}}
\bigoplus_{i=1}^2 \calh^{D_{\infty} \rtimes_a \bbZ}_n\left(
\underline{E}(D_{\infty}\rtimes_a \bbZ) \to \EGF{D_{\infty}\rtimes_a \bbZ}{\calvcyc};
\bfL^{\langle - \infty \rangle}(\bbZ?)\right)\Bigg)
\\
\xrightarrow{\cong}
\calh^G\left(\underline{E}G \to \EGF{G}{\calvcyc};\bfL^{\langle -\infty \rangle}(\bbZ?)\right),
\end{multline*}
where the term $\calh^{D_{\infty} \rtimes_a \bbZ}_n\left(
\underline{E}(D_{\infty}\rtimes_a \bbZ) \to
\EGF{D_{\infty}\rtimes_a \bbZ}{\calvcyc}; \bfL^{\langle - \infty
\rangle}(\bbZ?)\right)$ is analyzed in \eqref{relative term for D
rtimes Z};

\item \label{the: L-theory of hei rtimes Z/4: independent of decoration}
The canonical map $L^{\langle - \infty \rangle}(\bbZ G) \xrightarrow{\cong}
L^{\epsilon}(\bbZ G)$ is bijective for all decorations $\epsilon = p,h,s$.
 \end{enumerate}
\end{theorem}
\begin{proof}
~\ref{the: L-theory of hei rtimes Z/4: underline(E)}\qua This
follows from Theorem~\ref{the: $L$-theory}~\ref{the: $L$-theory: computation for underline E}
since  the groups
$\widetilde{L}_n^p(\bbZ/2) = \widetilde{L}_n^{\langle -\infty \rangle}(\bbZ/2)$ and
$\widetilde{L}_n^p(\bbZ/4) = \widetilde{L}_n^{\langle -\infty \rangle}(\bbZ/4)$ are torsionfree
\cite[Theorem 1]{Bak(1976)}.

\ref{the: L-theory of hei rtimes Z/4: splitting for fin and relative term}\qua
The Farrell--Jones Conjecture for algebraic $L$--theory with coefficients in $R = \bbZ$ is true
for $G = \Hei \rtimes \bbZ/4$ and $p^{-1}(V)$ for $V \subseteq Q$ virtually cyclic
since $G$ is a discrete cocompact subgroup of the virtually connected Lie group
$\Hei(\bbR) \rtimes_a \bbZ/4$ (see \cite{Farrell-Jones(1993a)}).
Since for a virtually cyclic group $V$ we have
$K_n(\bbZ V) = 0$ for $n \le -2$ \cite{Farrell-Jones(1995)},
we can apply
Theorem~\ref{the: $L$-theory}~\ref{the: $L$-theory: splitting due to Bartels}.
 
\ref{the: L-theory of hei rtimes Z/4: relative term}\qua
This follows from Lemma~\ref{lem: Z^n to Q to F and condition and calvcyc_f},
the lists \eqref{list of the V-s} and \eqref{list of the p^(-1)(V)-s} and
the isomorphism \eqref{relative term for D_(infty) times bbZ}.

\ref{the: L-theory of hei rtimes Z/4: independent of decoration}\qua
Because of the Rothenberg sequences it suffices to show that the Tate cohomology groups
$\widehat{H}^n(\bbZ/2,\Wh_q(G))$ vanish for $q \le 1$ and $n \in \bbZ$.
If $q \le -1 $, then $\Wh_q(G) = 0$, and,
if $q = 0,1$, then $\Wh_q(G) =  N\!K_q(\bbZ [\bbZ/4]) \bigoplus N\!K_q(\bbZ [\bbZ/4])$
by Corollary~\ref{cor: algebraic K-theory of G = Hei rtimes Z/4}. One easily checks
that the involution on $\Wh_q(G)$  corresponds under this identification to the involution
on $N\!K_q(\bbZ [\bbZ/4]) \bigoplus N\!K_q(\bbZ [\bbZ/4])$ which sends $(x_1,x_2)$ to
$(x_2,\tau(x_1))$ for  $\tau \colon N\!K_q(\bbZ [\bbZ/4]) \to N\!K_q(\bbZ [\bbZ/4])$
the involution on the Nil-Term. Hence the $\bbZ[\bbZ/2]$--module  $\Wh_q(G)$ is isomorphic
to the $\bbZ[\bbZ/2]$--module $\bbZ[\bbZ/2] \otimes_{\bbZ}  N\!K_q(\bbZ [\bbZ/4])$,
which is obtained from the $\bbZ$--module  $N\!K_q(\bbZ [\bbZ/4])$ by induction with the
inclusion of the trivial group into $\bbZ/2$. This implies
$\widehat{H}^n(\bbZ/2,\Wh_q(G)) = 0$  for $q \le 1$ and $n \in \bbZ$.
\end{proof}

\begin{remark} \label{rem: real K-theory and L-theory}
\em If one inverts $2$, then the computation for $L_n(\bbZ[\Hei \rtimes \bbZ])$
simplifies drastically as explained in the introduction because of
Lemma~\ref{lem: general simplification of calvcyc}.  In general this example shows how complicated
it is to deal with the infinite virtually cyclic subgroups which admit an epimorphism to $D_{\infty}$
and the resulting $\UNil$--terms.
\em
\end{remark}

\section{Group homology}
\label{sec: Group homology}

Finally we explain what the methods above give for the group homology

\begin{theorem} \label{the: group homology for the general case}
Let $G$ be the group appearing in \eqref{extension of groups} and assume
that conditions (M), (NM), and (T)  are satisfied. We then obtain
a long exact Mayer--Vietoris sequence
\begin{multline*}
\ldots \to H_{n+1}(G\backslash\underline{E}G) \xrightarrow{\partial_{n+1}}\\
\bigoplus_{i \in I} H_n(p^{-1}(M_i))
\xrightarrow{\left(\bigoplus_{i \in I} H_n(l_i)\right) \bigoplus
\left(\bigoplus_{i \in I} H_n(p^{-1}(M_i)\backslash s_i)\right)}
\\
H_n(G) \bigoplus \left(\bigoplus_{i \in I} H_n(p^{-1}(M_i)\backslash\underline{E}p^{-1}(M_i))\right)
\\
\xrightarrow{H_n(G\backslash s) \bigoplus \left(\bigoplus_{i \in I} H_n(d_i)\right)}
H_n(G\backslash\underline{E}G) \xrightarrow{\partial_n}
\\
\bigoplus_{i \in I} H_{n-1}( p^{-1}(M_i))
\xrightarrow{\left(\bigoplus_{i \in I} H_{n-1}(l_i)\right) \bigoplus
\left(\bigoplus_{i \in I} H_{n-1}(p^{-1}(M_i)\backslash s_i)\right)}  \ldots
\end{multline*}
from the pushout \eqref{pushout for G backslash underline E G} 
where $l_i \colon p^{-1}(M_i) \to G$ is the inclusion,
$s_i \colon Ep^{-1}(M_i) \to \underline{E}p^{-1}(M_i)$,
$s \colon EG \to \underline{E}G$ are the obvious equivariant maps and
$$d_i \colon p^{-1}(M_i)\backslash\underline{E}p^{-1}(M_i) \to G\backslash\underline{E}G$$ is
the map induced by the $l_i$--equivariant map $\underline{E}p^{-1}(M_i) \to \underline{E}G$.
\end{theorem}

\begin{remark} \label{rem: finite dimensional underline E G} \em
There are often finite-dimensional models for $\underline{E}G$ as discussed
in  \cite{Kropholler-Mislin(1998)}, \cite{Lueck(2000a)}. If or instance, there is a $k$--dimensional model
for $BK$ and a $m$--dimensional model for $\underline{E}Q$ and $d$ is a positive integer
such that the order of any finite subgroup of $Q$ divides $d$, then there is a
$(dk+n)$--dimensional model for $\underline{E}G$  \cite[Theorem 3.1]{Lueck(2000a)}.
If $Q$ is an extension $0 \to \bbZ^n \to Q \to F \to 1$ for a finite group $F$
and there is a $k$--dimensional model for $BK$, then there is a
$(|F|\cdot k+n)$--dimensional model for $\underline{E}G$.

Suppose that there is a $N$--dimensional model for $\underline{E}G$.
Then there is also a $N$--dimensional model for $\underline{E}p^{-1}(M_i)$
for each $i \in I$ and under the assumptions of
Theorem~\ref{the: group homology} we obtain for $n \ge N+1$ an isomorphism
$$\bigoplus_{i \in I} H_n(l_i) \colon \bigoplus_{i \in I} H_n(p^{-1}(M_i))
\xrightarrow{\cong} H_n(G).$$
\em
\end{remark}

Next we compute the group homology $H_*(\Hei \rtimes \bbZ/4)$.
We start with the computation of $H_n(\Hei)$.
The Atiyah-Hirzebruch spectral sequence associated to the central extension
$1 \to \bbZ \xrightarrow{i'}\Hei \xrightarrow{p'} \bbZ^2 \to 0$ yields the isomorphism
$$H_2(\bbZ^2)  \xrightarrow{\cong} H_3(\Hei)$$
and the long exact sequence
\begin{multline*}
0 \to H_1(\bbZ^2) \to H_2(\Hei) \xrightarrow{H_2(p')} H_2(\bbZ^2) \to H_0(\bbZ^2) = H_1(S^1) = H_1(\bbZ)
\\
\xrightarrow{H_1(i')}  H_1(\Hei) \xrightarrow{H_1(p')} H_1(\bbZ^2) \to 0.
\end{multline*}
Since $z \in \Hei$ is a commutator, namely $[u,v]$, the map $H_1(i') \colon H_1(\bbZ) \to
H_1(\Hei)$ is trivial. This implies:

\begin{lemma} \label{lem: homology of BHei}
There are natural isomorphisms
$$\begin{array}{rclcl}
H_1(p') \colon H_1(\Hei) & \xrightarrow{\cong} & H_1(\bbZ^2); & &
\\
H_1(\bbZ^2) &  \xrightarrow{\cong} & H_2(\Hei); & &
\\
H_2(\bbZ^2) & \xrightarrow{\cong} & H_3(\Hei);
\\
H_n(\Hei) & = & 0 & & \text{ for n } \ge 4.
\end{array}$$
\end{lemma}

Next we analyze the Atiyah-Hirzebruch spectral sequence associated
to the split extension $1 \to \Hei \xrightarrow{k} G := \Hei
\rtimes \bbZ/4 \xrightarrow{\pr} \bbZ/4 \to 1$. The isomorphisms
above appearing in the computation of the homology of $\Hei$ are
compatible with the $\bbZ/4$--actions. Thus we get
$$\begin{array}{lclcl}
H_p(\bbZ/4;H_q(\Hei)) & = & H_p(\bbZ/4;H_q(\bbZ^2)) = \bbZ/2 & & \text{ for }
q= 1,2, p \ge 0, p \text{ even};
\\
H_p(\bbZ/4;H_q(\Hei)) & = & H_p(\bbZ/4;H_q(\bbZ^2)) = 0 & & \text{ for }
q= 1,2, p \ge 0, p \text{ odd};
\\
H_p(\bbZ/4;H_q(\Hei)) & = & H_p(\bbZ/4) & & \text{ for }
q= 0,3;
\\
H_p(\bbZ/4;H_q(\Hei)) & = & = 0 & & \text{ for }
q \ge 4.
\end{array}$$
Hence the $E^2$--term looks like:
$$\begin{array}{ccccccc}
\bbZ & \bbZ/4 & 0 & \bbZ/4 & 0 & \bbZ/4 & 0
\\
\bbZ/2  & 0  & \bbZ/2 & 0 & \bbZ/2 & 0 & \bbZ/2
\\
\bbZ/2  & 0  & \bbZ/2 & 0 & \bbZ/2 & 0 & \bbZ/2
\\
\bbZ & \bbZ/4 & 0 & \bbZ/4 & 0 & \bbZ/4 & 0
\end{array}$$
Using the model for $\underline{E}G$ of
Lemma~\ref{G backslash underline E G is S^3} we see that the map
$B\!\Hei \to G\backslash \underline{E}G$ can be identified with the quotient map
$B\!\Hei \to \bbZ/4\backslash B\!\Hei$ of an orientation preserving smooth $\bbZ/4$--action
on the closed orientable $3$--manifold $B\!\Hei$, where the quotient is again a closed
orientable $3$--manifold and the action has at least one free orbit. Since we can compute
the degree of a map by counting preimages of a regular value, the degree
must be $\pm 4$. Recall that $G \to \bbZ/4$ is split surjective.
These remarks imply together with the spectral sequence above

\begin{lemma} \label{H_3(Hei) to H_3(G)}
The composition $H_3(B\!\Hei) \xrightarrow{H_3(Bk)} H_3(BG)
\xrightarrow{H_3(G\backslash s)} H_3(G\backslash\underline{E}G)$
is an injective map of infinite cyclic subgroups whose cokernel has order four.
The map $H_3(\Hei) \to H_3(G)$ is injective and the order of the cokernel of the induced map
$H_3(G)/\tors(H_3(G)) \to H_3(G\backslash \underline{E}G)$ divides four;

Moreover, there are the following possibilities
\begin{enumerate}
\item The differential $d^2_{2,1} \colon E^2_{2,1} \cong \bbZ/2 \to E^2_{0,2} \cong
  \bbZ/2$ is trivial. Then $H_2(G)$ is $\bbZ/2$. Moreover, either
the group $H_3(G)$ is $\bbZ \times \bbZ/4$ and the
induced map $H_3(\Hei) \to H_3(G)/\tors(H_3(G))$ is
an injective homomorphism of infinite cyclic groups
whose cokernel has order two, or the group $H_3(G)$ is $\bbZ \times \bbZ/2 \times \bbZ/4$
and the induced map
$H_3(\Hei) \to H_3(G)/\tors(H_3(G))$ is an isomorphism of infinite cyclic groups.

\item The differential $d^2_{2,1} \colon E^2_{0,2} \cong \bbZ/2 \to E^2_{2,1} \cong
  \bbZ/2$ is non-trivial. Then $H_2(G)$ is $0$.

\end{enumerate}
\end{lemma}

It is not obvious how to compute the homology groups $H_n(G)$ for
$G = \Hei\rtimes\bbZ/4)$ from the Atiyah-Hirzebruch spectral
sequence. Let us try Theorem~\ref{the: group homology for the
general case}. It yields the long exact Mayer Vietoris sequence
\begin{multline*}
\ldots \to H_{n+1}(G\backslash\underline{E}G) \xrightarrow{\partial_{n+1}}
\bigoplus_{i = 0}^2 H_n(p^{-1}(M_i))
\\
\xrightarrow{\left(\bigoplus_{i = 0}^2 H_n(l_i)\right) \bigoplus
\left(\bigoplus_{i = 0}^2 H_n(p^{-1}(M_i)\backslash s_i)\right)}
H_n(G) \bigoplus \left(\bigoplus_{i = 0}^2 H_n(p^{-1}(M_i)\backslash\underline{E}p^{-1}(M_i))\right)
\\
\xrightarrow{H_n(G\backslash s) \bigoplus \left(\bigoplus_{i = 0}^2 H_n(d_i)\right)}
H_n(G\backslash\underline{E}G) \xrightarrow{\partial_n}
\\
\bigoplus_{i = 0}^2 H_{n-1}( p^{-1}(M_i))
\xrightarrow{\left(\bigoplus_{i = 0}^2 H_{n-1}(l_i)\right) \bigoplus
\left(\bigoplus_{i = 0}^2 H_{n-1}(p^{-1}(M_i)\backslash s_i)\right)}  \ldots
\end{multline*}
where the maximal finite subgroups $M_0$, $M_1$ and $M_2$ of $Q$ have been introduced
in Lemma~\ref{lem: maximal finite subgroups of Q}.
The map $s_i \colon p^{-1}(M_i)\backslash Ep^{-1}(M_i) \to p^{-1}(M_i)\backslash
\underline{E}p^{-1}(M_i)$ can be identified with
\begin{eqnarray*}
s_0\colon B\langle t,z\rangle = B\langle t \rangle \times B\langle z \rangle
& \xrightarrow{\pr} & B\langle z \rangle;
\\
s_1 \colon B\langle ut \rangle & \xrightarrow{\id} & B\langle ut \rangle;
\\
s_2\colon B\langle ut^2,z \rangle = B\langle ut^2\rangle \times B\langle z \rangle
& \xrightarrow{\pr} & B\langle z \rangle.
\end{eqnarray*}
Hence we obtain the exact sequence
\begin{multline*}
\ldots \to H_{n+1}(G)   \bigoplus H_{n+1}(\langle z \rangle ) \bigoplus H_{n+1}(\langle ut  \rangle) \bigoplus
H_{n+1}(\langle z \rangle)
\\
\xrightarrow{H_{n+1}(G\backslash s) \bigoplus_{i=0}^2 H_{n+1}(d_i)}
H_{n+1}(G\backslash\underline{E}G)
\\
\xrightarrow{\partial_{n+1}^{\prime}}
H_n(\langle t \rangle \times \langle z \rangle) \bigoplus H_n(\langle ut  \rangle)
\bigoplus H_n(\langle ut^2  \rangle \times \langle z \rangle)
\\
\xrightarrow{\left(H_n(\incl_0^{\prime}) \bigoplus  H_n(\incl_1^{\prime})
\bigoplus H_n(\incl_2^{\prime})\right)
\bigoplus H_n(\pr) \bigoplus \id \bigoplus H_n(\pr)}
\\
H_n(G) \bigoplus H_n(\langle z \rangle ) \bigoplus H_n(\langle ut  \rangle) \bigoplus
H_n(\langle z \rangle)\\
\xrightarrow{H_n(G\backslash s) \bigoplus_{i=0}^2 H_n(d_i)} H_n(G\backslash\underline{E}G)
\xrightarrow{\partial_n'} \ldots
\end{multline*}
where
\begin{eqnarray*}
\incl_0^{\prime} \colon \langle t,z \rangle = \langle t \rangle \times \langle z \rangle
& \to & G
\\
\incl_1^{\prime} \colon \langle ut\rangle & \to & G
\\
\incl_2^{\prime} \colon \langle ut, z\rangle = \langle ut \rangle \times \langle z \rangle
& \to & G
\end{eqnarray*}
are the inclusions. This yields the exact sequence:
\begin{multline}
\ldots \to H_{n+1}(G)
\xrightarrow{H_{n+1}(G\backslash s)}
H_{n+1}(G\backslash\underline{E}G)
\\
\xrightarrow{\partial_n''}
\widetilde{H}_n(\langle t \rangle)  \bigoplus
\widetilde{H}_{n-1}(\langle t \rangle) \bigoplus \widetilde{H}_n(\langle ut^2\rangle)
\bigoplus \widetilde{H}_{n-1}(\langle ut^2\rangle)
\\\to
H_n(G)  \xrightarrow{H_n(G\backslash s)}
H_{n}(G\backslash\underline{E}G) \xrightarrow{\partial_n''} \ldots
\label{exact sequence for H_*(G)}
\end{multline}
Recall that $G\backslash\underline{E}G$ is $S^3$.
We conclude from  Lemma~\ref{H_3(Hei) to H_3(G)} that the order of the cokernel of the map
$H_3(G\backslash s) \colon H_3(BG) \to H_3(G\backslash \underline{E}G)$ divides four.
Since the order of
$\widetilde{H}_2(\langle t \rangle)  \bigoplus
\widetilde{H}_1(\langle t \rangle) \bigoplus \widetilde{H}_2(\langle e_1t^2\rangle)
\bigoplus \widetilde{H}_2(\langle e_1t^2\rangle)$ is eight, the long exact sequence above
implies that the group $H_2(G)$ is
different from zero and that the group $H_3(G)$ is isomorphic to $\bbZ \times \bbZ/2
\times \bbZ/4$.
Now Lemma~\ref{H_3(Hei) to H_3(G)} and the long exact sequence
\eqref{exact sequence for H_*(G)} above imply

\begin{theorem} \label{the: group homology}
For $G = \Hei \rtimes \bbZ/4$ we have isomorphisms
\begin{eqnarray*}
H_n(G) & = & \bbZ/2 \times \bbZ/4  \text{ for } n \ge 1, n \not= 2,3;
\\
H_2(G) & = & \bbZ/2;
\\
H_3(G) & = & \bbZ \times \bbZ/2 \times \bbZ/4.
\end{eqnarray*}
The map $H_3(\Hei) \to H_3(G)/\tors(H_3(G))$ is an isomorphism.
\end{theorem}

One can compute the group cohomology analogously or derive it from the homology
by the universal coefficient theorem.

\section{Survey over other extensions}
\label{sec: Survey over other extensions}

There are other prominent extensions of $\Hei$ which can be treated analogously to the case
$\Hei \rtimes \bbZ/4$. We give a brief summary of the topological $K$--theory and the algebraic $K$--theory below.
In all cases $G\backslash \underline{E}G$ is $S^3$.

\subsection{Order six symmetry}
\label{sec: Order six symmetry}

Consider the following automorphism $\omega \colon \Hei \to \Hei$ of order $6$ which sends
$u$ to $v$, $v$ to $u^{-1}v$, and $z$ to $z$. 

\begin{theorem}
\label{the: K_1(C^*_r(G)) in C^*-terms for G = Hei rtimes Z/6}
For the group
\begin{multline*}
G = \Hei \rtimes \bbZ/6 = \langle u,v,z,t \mid [u,v] = z, t^6 = 1, [u,z] = [v,z] =
[t,z] = 1, \\
tut^{-1} = v, tvt^{-1} = u^{-1}v \rangle
\end{multline*}
there is a short exact sequence
$$
0 \to\widetilde{R}_{\bbC}(\langle t \rangle)  \to
K_1(C^*_r(G)) \to
\widetilde{K}_1(S^3) \to 0
$$
and an isomorphism
$$
R_{\bbC}(\langle t \rangle)  \xrightarrow{\cong} 
K_0(C^*_r(G)).
$$
There are isomorphisms
\begin{eqnarray*}
\Wh_n(G) & \cong & \left\{\begin{array}{lll}
N\!K_1(\bbZ [\bbZ/6]) \bigoplus N\!K_1(\bbZ [\bbZ/6])  & & \text { for  } n = 1;
\\
\Wh_{-1}(\bbZ/6) \cong \bbZ & & \text{ for } n = -1,0;
\\
0 & & \text{ for } n \le -2.
\end{array}\right.
\end{eqnarray*}

\end{theorem}


\subsection{Order three symmetry}
\label{sec: Order three symmetry}

Next we deal with the $\bbZ/3$--action on $\Hei$ given by $\omega^2$,
where $\omega$ is the automorphism of order six investigated in
Subsection~\ref{sec: Order six symmetry}. 

\begin{theorem}
\label{the: K_1(C^*_r(G)) in C^*-terms for G = Hei rtimes Z/3}
For the group 
\begin{multline*}
G = \Hei \rtimes \bbZ/3 = \langle u,v,z,t \mid [u,v] = z, t^3 = 1, [u,z] = [v,z] =
[t,z] = 1, \\
tut^{-1} = u^{-1}v, tvt^{-1} = u^{-1} z^{-1}\rangle
\end{multline*}
there is a short exact sequence
$$
0 \to \widetilde{R}_{\bbC}(\langle t \rangle) \to K_1(C^*_r(G)) \to \widetilde{K}_1(S^3) \to 0
$$
and an isomorphism
$$R_{\bbC}(\langle t \rangle) \xrightarrow{\cong} K_0(C^*_r(G)).$$
We have $\Wh_n(G) = 0$ for $n \le 2$.

The $L$--groups $L_n{\epsilon}(\bbZ G)$ are independent of the choice of decoration
$\epsilon = -\infty,$ $p ,h,s$ and the reduced ones fit into a short split exact sequence
\begin{multline*}
0 \to \widetilde{L}_n^{\langle - \infty \rangle}(\bbZ\langle t\rangle) \bigoplus
 \widetilde{L}_{n-1}^{\langle - \infty \rangle}(\bbZ\langle t\rangle) \to
\widetilde{L}_n^{\langle - \infty \rangle}(\bbZ G) \to L_{n-3}^{\langle - \infty \rangle}(\bbZ) \to 0.
\end{multline*}
\end{theorem}


\subsection{Order two symmetry}
\label{sec: Order two symmetry}

Next we deal with the $\bbZ/2$--action on $\Hei$ given by $u \mapsto u^{-1}$,
$v \mapsto v^{-1}$ and $z \mapsto z$. This is the square of the automorphism of order four used in
the $\bbZ/4$--case. 

\begin{theorem}
\label{the: K_1(C^*_r(G)) in C^*-terms for G = Hei rtimes Z/2}
For the group 
\begin{multline*}
G = \Hei \rtimes \bbZ/2 = \langle u,v,z,t \mid [u,v] = z, t^2 = 1, [u,z] = [v,z] =
[t,z] = 1, \\
tut^{-1} = u^{-1}, tvt^{-1} = v^{-1} \rangle
\end{multline*}
there is a short exact sequence
$$
0 \to \oplus_{i=0}^2 \widetilde{R}_{\bbC}(M_i) \to K_1(C^*_r(G)) \to \widetilde{K}_1(S^3) \to 0
$$
and an isomorphism
$$
0 \to K_0(\pt) \bigoplus \bigoplus_{i=0}^2 \widetilde{R}_{\bbC}(M_i)
\xrightarrow{\cong }
K_0(C^*_r(G)),
$$
where
\begin{eqnarray*}
M_0 & = & \langle t \rangle ;\\
M_1 & = &  \langle ut \rangle;\\
M_2 & = &  \langle vt \rangle.
\end{eqnarray*}
We have $\Wh_n(G) = 0$ for $n \le 2$.
\end{theorem}

\end{document}